\documentclass[review]{elsarticle}

\usepackage{amsmath}	
\usepackage{amssymb}
\usepackage{algorithm}
\usepackage{algorithmic}
\usepackage{color}
\newcommand{\bs}[1]{\boldsymbol{#1}}

\newcommand{\figref}[1]{{Fig. \ref{#1}}}

\begin{document}
\begin{frontmatter}

\title{A Burton-Miller-type boundary element method based on a hybrid integral representation and its application to cavity scattering}

\author[1]{Riku Toshimitsu}
\author[1]{Hiroshi Isakari\corref{cor1}}
\ead{isakari@sd.keio.ac.jp}

\cortext[cor1]{Corresponding author}

\address[1]{{Faculty of Science and Technology, Keio University},
{3-14-1, Hiyoshi, Kohoku-ku},
{Yokohama, Kanagawa},
{223-8522},
{Japan}}

\begin{abstract}
This study builds on a recent paper by Lai et al [{\it Appl. Comput. Harmon. Anal.}, 2018] in which a novel boundary integral formulation is presented for scalar wave scattering analysis in two-dimensional layered and half-spaces. The seminal paper proposes a hybrid integral representation that combines the Sommerfeld integral and layer potential to efficiently deal with the boundaries of infinite length. In this work, we modify the integral formulation to eliminate the fictitious eigenvalues by employing Burton-Miller's approach. We also discuss reasonable parameter settings for the hybrid integral equation to ensure efficient and accurate numerical solutions. Furthermore, we extend the modified formulation for the scattering from a cavity in a half-space whose boundary is locally perturbed. To address the cavity scattering, we introduce a virtual boundary enclosing the cavity and couple the integral equation on it with the hybrid equation. The effectiveness of the proposed method is demonstrated through numerical examples.
\end{abstract}


\end{frontmatter}


\section{Introduction}
Sound propagation is a topic of interest in various aspects of daily life and industry. For instance, any products ranging from daily necessities such as aeroplanes, trains, cars, and household appliances to huge industrial machines must carefully be designed to reduce noise. Not just reducing the noise, we sometimes need more sophisticated sound control in designing for example a concert hall to provide a better experience. As a tool for such designs, numerical methods for solving the acoustic scattering problem, a mathematical model of sound propagation, have long been studied by many researchers and engineers. As widely accepted numerical methods for the problems, we can mention the finite-difference time-domain (FDTD) method, finite element method (FEM) and boundary element method (BEM). Among them, the BEM would be the best choice when an infinite space is the target. This is partly because the BEM requires boundary-only discretisation whereas the FDTD and FEM need to discretise the target space. In other words, the BEM can rigorously address the scattering in infinite space as long as the involved boundary is finite-sized.

If an infinitely large boundary is involved, on the other hand, it may be difficult to use a naive BEM to solve the scattering problem. To address such a problem, e.g. acoustic scattering in a semi-infinite domain, a classical approach by using the Sommerfeld integral is sometimes useful~\cite{sommerfeld1909uber,vanderpol1935theory}. One may transform the layer potential defined over the infinite measure into the Sommerfeld integral by replacing the kernel function with its Fourier transform along the infinite boundary (and accordingly the density function is also transformed into a Fourier space). In some cases, the resulting integral can accurately be approximated by an integral over a finite domain. This approach fails, however, when for example either the point source or the scatterer is close to the infinite boundary. In such cases, even after truncation, an impractically wide integral range is needed for the Sommerfeld integral to be accurately evaluated. Other possible methods for scattering in semi-infinite or layered domains include the method of complex images~\cite{ochmann2004complex} and some special integral formulations with appropriate Green's functions~\cite{perez-arancibia2014highorder}. Each method has, however, its own limitations. Further researches are still ongoing to establish an efficient numerical method for these problems. Among such efforts, a recently proposed integral formulation~\cite{lai2018new} seems to be promising. The novel formulation is based on the hybrid integral representation that utilises both the layer potential in physical space and the Sommerfeld integral in the Fourier space. By ``appropriately'' truncating both the integral domains, all the relevant integrals can be defined over a finite domain of moderate size without sacrificing accuracy for the scattered field even when neither the layer potential nor the Sommerfeld integral can. It is unfortunate, however, that the original hybrid integral representation may suffer from the fictitious eigenvalue problem. One of the main goals of this paper is thus to propose a well-defined integral representation that is free from the problem. We here present a new hybrid integral representation combined with the standard Burton-Miller approach~\cite{burton1971application} to this end. We also investigate the truncation parameters for the hybrid representation. In order for the method to efficiently and accurately work, both the finite integral domains for the layer potential and the Sommerfeld integral should carefully be set. In the present work, we sweep the parameters to find reasonable settings for some benchmark problems.

Another main focus of this study is to extend the hybrid approach to address so-called {\it cavity scattering}, which deals with wave scattering in a semi-infinite space whose boundary is locally perturbed in a concave shape. The cavity scattering is of engineering interest in a wide range of fields. The possible applications may include laser sensing~\cite{bao2014optimal,bao2015radar} in electromagnetic scattering, and it is also expected to play an important role in acoustic sensing and/or advanced acoustic devices that utilise sound reflection and resonance in open space. It is, however, due to the concaved cavity shape, difficult to adopt the standard method of images. So far, several boundary integral formulations for the problem have been proposed; Ammari et al~\cite{ammari2000integral}, for example, introduced a virtual boundary to cover the existing cavity, and coupled a boundary integral formulation in the cavity and the Sommerfeld integral, whose unique solvability is mathematically justified. Since it is basically based on the Sommerfeld integral, it may need a large truncated integral to ensure an accurate solution. Lai et al~\cite{lai2014fast} also tackled the problem, in which the integral equation on the virtual boundary as well as the method of image is utilised. The formulation involves, however, somewhat technical insight from viewpoint of physics to have a well-posed integral equation. In this paper, we propose a simple boundary integral formulation for the cavity scattering that couples the standard boundary integral on the virtual boundary and the hybrid integrals. 

The rest of the paper is organised as follows: Section 2 discusses the acoustic scattering in half-space introducing the hybrid BEM combined with the Burton-Miller formulation and the parameter settings for the hybrid BEM with some validation results. Section 3 is devoted to cavity scattering analysis, where some configurations demonstrated in the numerical examples might work as open resonators. Section 4 concludes the paper. The Appendix summarises the hybrid integral representation for referential purposes. 

\section{A new BEM based on the hybrid integral representation of Burton-Miller type}
\subsection{Settings}\label{sec:setting2}
\begin{figure}[h]
  \centering
  \includegraphics[scale=0.35]{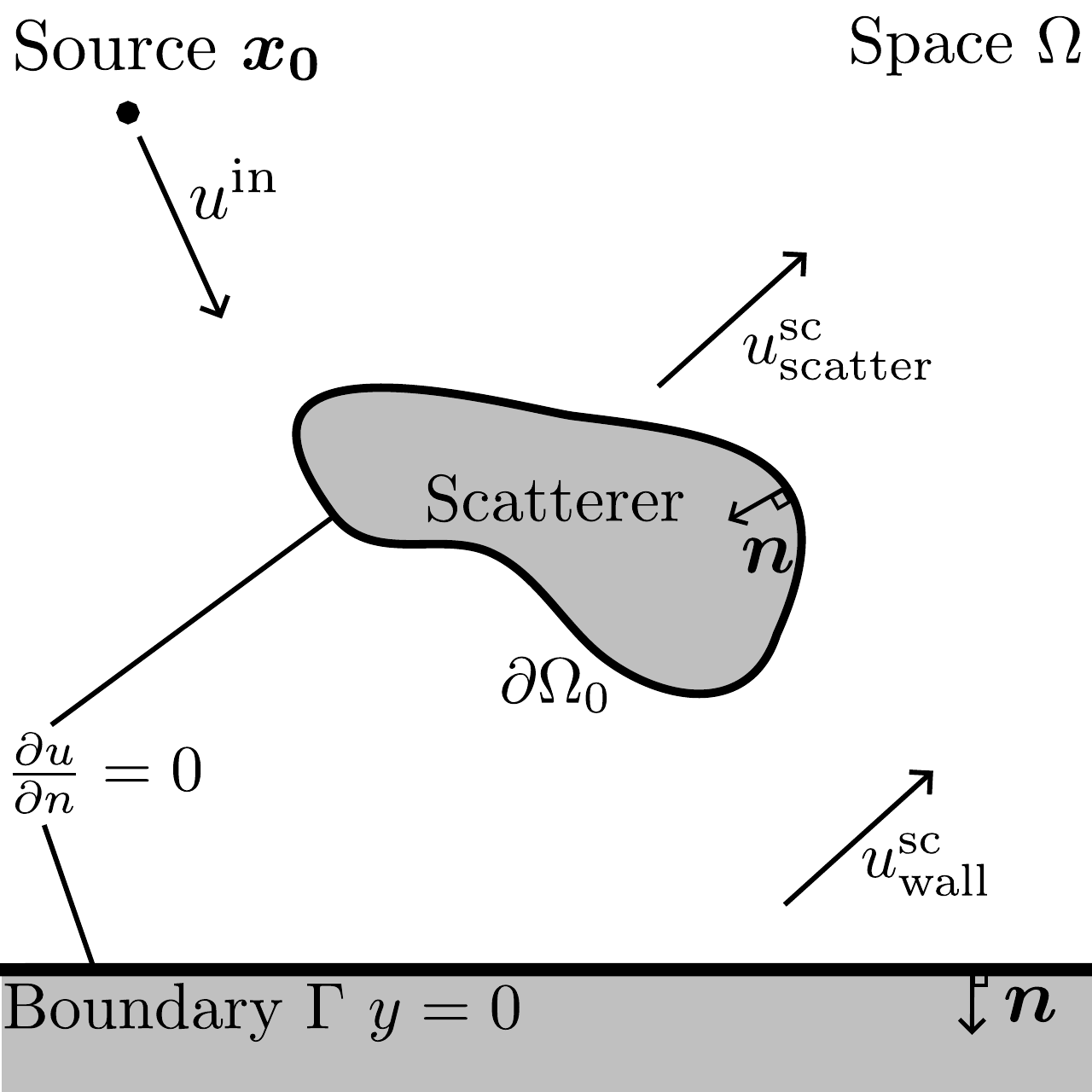}
  \caption{Illustrative sketch for the acoustic scattering problem in 2D half-space.}\label{fig:setting2}
\end{figure}

Let us here consider the wave scattering in a two-dimensional half-space (\figref{fig:setting2}), in which an infinitely long wall $\Gamma:=\{\bs{x}=(x, y) \mid y=0\}$ and a bounded scatterer $\Omega_0$ are arranged. Here, both the wall and scatterer are assumed to be acoustically rigid. We also put a sound source at $\bs{x_0}$ in $\Omega:=\mathbb{R}^2_+\setminus\overline{\Omega_0}$, where $\mathbb{R}^2_+:=\{\bs{x}\mid y>0\}$ is the upper half-space. The incident field induced by the source vibrating with the angular frequency $\omega$ is denoted as $u^\mathrm{in}(\bs{x})=\mathrm{i} H_0^{(1)}(k|\bs{x}-\bs{x_0}|)/4$, where $k:=\omega/c$ is the wave number, and $H_n^{(1)}$ is the Hankel function of the first kind and $n^\mathrm{th}$ order, and $\mathrm{i}$ is the imaginary unit. We assume that a compressive and inviscid fluid with phase velocity $c$ is filled in $\Omega$. We are interested in finding the scattered field $u^\mathrm{sc}_\mathrm{wall}+u^\mathrm{sc}_\mathrm{scatterer}$, where $u^\mathrm{sc}_\mathrm{wall}$ and $u^\mathrm{sc}_\mathrm{scatterer}$ are respectively the scattered wave by the wall and scatterer. 
With these settings, the total field 
\begin{align}
 u(\bs{x})=u^{\mathrm{in}}(\bs{x})+u^{\mathrm{sc}}_{\mathrm{wall}}(\bs{x})+u^{\mathrm{sc}}_{\mathrm{scatterer}}(\bs{x})
\label{eq:TOTALFIELD}
\end{align}
is governed by the following Neumann boundary value problem (BVP) of the Helmholtz equation in two dimension:
\begin{align}
 \nabla^2 u(\bs{x}) +k^2 u(\bs{x}) + \delta(\bs{x}-\bs{x_0}) = 0 ~~\bs{x}\in \Omega, \label{eq:Helmholtz}\\
 \frac{\partial u(\bs{x})}{\partial \bs{n}(\bs{x})} =0  ~~\bs{x} \in \Gamma\cup\partial\Omega_0, \label{eq:Neumann}\\
 \mathrm{Radiation~condition~for~}~u(\bs{x})-u^{\mathrm{in}}(\bs{x})~~\mathrm{as}~~|\bs{x}|\to \infty, \label{eq:radiation}
\end{align}
where $\delta$ is the Dirac delta, and the unit normal $\bs{n}$ is directed from $\Omega$. Note that, with the solution of the BVP, the total sound pressure in time domain is recovered as $\Re[u(\bs{x})e^{-\mathrm{i}\omega t}]$ at time $t$. 

\subsection{The Burton-Miller integral formulation}
We here formulate a boundary element method for the scattering problem defined in Section \ref{sec:setting2}. In accordance with the discussions in Appendix A, we use the following integral representation for the scattered wave from the infinite boundary $\Gamma$: 
\begin{align}
  u^{\mathrm{sc}}_{\mathrm{wall}}(\bs{x})  =S_{\Gamma_0}[\sigma_W](\bs{x}) +F_{I_0}[ \widehat{\xi_W}](\bs{x}), 
\label{eq:SSSO}
\end{align}
with the integral operators $S_{\Gamma_0}$ and $F_{I_0}$ respectively defined as \eqref{eq:sec2_4_1} and \eqref{eq:sec2_4_2}. Also, the truncated density function $\sigma_W$ is accordingly defined as \eqref{eq:sgmW} with the window function \eqref{eq:window}, and $\widehat{\xi_W}$ is the unknown density function for the Sommerfeld integral. The finite intervals $\Gamma_0=\{\bs{x}= \mid -M_0 < x < M_0,~y=0 \}$ and $I_0=(-N_0,~N_0)$ are also introduced as in Appendix A, with $0< M_0, N_0<\infty$. For the scattering from the rigid inclusion $\Omega_0$, we use the following standard ansatz:
\begin{align}
u^{\mathrm{sc}}_{\mathrm{scatterer}}(\bs{x})=S_{\partial \Omega_0}[\sigma](\bs{x})+\beta D_{\partial \Omega_0}[\sigma](\bs{x}), 
\label{eq:bmintegral}
\end{align}
with the unknown density function $\sigma$, where $S_{\partial\Omega_0}$ is the operator for the single layer potential \eqref{eq:sec2_4_1}, and $D_{\partial \Omega_0}$ is that for the double layer one defined as
\begin{align}
 D_{\partial \Omega_0}[\sigma](\boldsymbol{x}):= \int_{\partial \Omega_0}  \frac{\partial G_k(\boldsymbol{x},\boldsymbol{x'})}{\partial \bs{n}(\bs{x'})} \sigma|_{\partial \Omega_0}(\bs{x'}) ds(\boldsymbol{x'}) \label{7}, 
\end{align}
where $G_k$ is the fundamental solution defined as \eqref{eq:fundamentalsol}.
In \eqref{eq:bmintegral}, $\beta\in\mathbb{C}$ is a parameter for the Burton-Miller method~\cite{burton1971application}. It should be noted that the original integral representation~\cite{lai2018new}, which may suffer from the fictitious eigenvalue problem as shall be seen later, corresponds to the case of $\beta=0$ in ours. It should also be noted that the Burton-Miller representation for \eqref{eq:SSSO} is not necessary because it is defined on an infinitely large and open boundary.

The remaining task is to find the density functions $\sigma$ on $\Gamma_0\cup\partial \Omega_0$ and $\widehat{\xi_W}$ on $I_0$. Corresponding to the local integral equation \eqref{eq:aug_int_eq}, we impose $\sigma$ to satisfy the following integral equation for $\bs{x}\in\Gamma_0$:
\begin{align}
  \frac{1}{2} \sigma|_{\Gamma_0}(\bs{x})+ D^*_{\Gamma_0}[\sigma](\bs{x})+D^*_{\partial \Omega_0}[\sigma](\bs{x})+\beta N_{\partial \Omega_0}[\sigma](\bs{x})=-\frac{\partial G_k(\bs{x},\bs{x_0})}{\partial \bs{n}(\bs{x})},\label{eq:H1BM_1}
\end{align}
where $D^*_{\Gamma_0}$ is the integral operater defined as \eqref{eq:dlayer}, and $N_{\nu}[\sigma]$ represents the following hyper singular integral operator:
\begin{equation}
  \begin{split}
      N_{\nu}[\sigma](\bs{x})&:=\frac{\partial D_{\nu}[\sigma](\bs{x})}{\partial \bs{n}(\bs{x})}=  \int_{\nu}  \hspace{-12pt}= \frac{\partial G_k(\boldsymbol{x},\boldsymbol{x'})}{\partial \bs{n}(\bs{x})\partial \bs{n}(\bs{x'})}\sigma (\boldsymbol{x'})ds(\boldsymbol{x'}),
  \end{split}
\end{equation}
where $\displaystyle\int \hspace{-12pt}=$ indicates the finite part of diverging integral. In order for the total field \eqref{eq:TOTALFIELD} with the integral representations \eqref{eq:SSSO} and \eqref{eq:bmintegral} to satisfy the boundary condition \eqref{eq:Neumann} on $\partial\Omega_0$, the density functions need to solve the following integral equation for $\bs{x}\in\partial\Omega_0$:
\begin{align}
  D^*_{\Gamma_0}[\sigma_W](\bs{x})+H_{I_0}[\widehat{\xi_W}](\bs{x})+\frac{1}{2} \sigma|_{\partial \Omega_0}(\bs{x})+D^*_{\partial \Omega_0}[\sigma](\bs{x}) +\beta N_{\partial \Omega_0}[\sigma](\bs{x})=-\frac{\partial G_k(\bs{x},\bs{x_0})}{\partial \bs{n}(\bs{x})},\label{eq:H1BM_2}
\end{align}
where $H_{I_0}[\widehat{\xi_W}]$ is related to the normal derivative of $F_{I_0}[\widehat{\xi_W}]$ and defined as
\begin{align}
  H_{I_0}[\widehat{\xi_W}](\bs{x}):=\frac{1}{4\pi}\int_{I_0}(i \lambda n_{x}-\sqrt{\lambda^2 -k^2}n_{y}) \frac{e^{- \sqrt{\lambda^2 - k^2 }y}}{\sqrt{\lambda^2 -k^2}}e^{i \lambda x} \widehat{\xi_W} (\lambda) d\lambda. \label{eq:H_I_0}
\end{align}
By the boundary condition \eqref{eq:Neumann} on $\Gamma$ (including $\Gamma_0$), we have the following integral equation:
\begin{equation}
 \begin{split}
  &\frac{1}{2} \widehat{\sigma_W}(\lambda,0)+\widehat{D^{*}_{\Gamma_0}[\sigma_W]}(\lambda,0) +\frac{1}{2} \widehat{\xi_W}   (\lambda)  \\ 
  &+ \widehat{D^{*}_{\partial \Omega_0}[\sigma]}(\lambda,0)
        +\beta \widehat{N_{\partial \Omega_0}[\sigma]}(\lambda,0)=-\widehat{\frac{\partial G}{\partial \bs{n}(\bs{x})}}(\lambda,0,\bs{y_0}), \label{eq:H1BM_3}
 \end{split}
 \end{equation}
where $\hat{\cdot}$ indicates the Fourier transform (see \eqref{eq:def_fourier} for its definition). The Fourier transforms of the layer potentials are evaluated by exploiting the convolutional structure of the potentials.

After discretisation by e.g. the collocation, one can solve the system of integral equations \eqref{eq:H1BM_1}, \eqref{eq:H1BM_2} and \eqref{eq:H1BM_3} simultaneously for the unknown density functions. Note that, while the integral equation \eqref{eq:aug_int_eq} for the case of $\Omega_0=\emptyset$ is decoupled from the equation in the Fourier space \eqref{eq:sol_global_equation}, \eqref{eq:H1BM_1} with $\Omega_0$ is not from \eqref{eq:H1BM_3}. 

\subsection{Validation}\label{sec:seido1}
In this subsection, we validate the proposed BEM with the Burton-Miller representation by checking if the method is free from the fictitious eigenvalue problem. To this end,  we compare the results from our method with those from the conventional method of image. As the accuracy indicator, we here use the following relative error:
\begin{align}
\frac{\Sigma_{j} |p_j-c_j|}{\Sigma_{j} |c_j|}, \label{eq:Relativeerror}
\end{align}
where $p_j$ and $c_j$ represent the value of $u$ at $j^\mathrm{th}$ interior point in $\Omega$ obtained by the proposed and conventional approaches, respectively. As a benchmark setting, we here consider the following: A circular rigid scatterer with the radius of $1.0$ centred at $(0.0,1.5)$ is allocated in $\mathbb{R}^2_+$, the point source is set at $\bs{x_0}=(1.0,3.0)$, oscillating with the angular frequency $\omega\in[0.01, 10.]$. With these settings, we compute $u$ at 10201 interior points distributed in a rectangle $\{\bs{x} \mid -4.0 \le x \le 4.0,~0.10 \le y \le 8.0\}$ on equal intervals. As for the numerical parameters, we set $a=2.0$ for the contour deformation \eqref{eq:contour_deformation}, $M_0=20$ and $N_0=30$ respectively for the intervals for layer potential on the wall and the Fourier integral. Note that we use sufficiently large integral intervals for this computation since we here are concerned only about the accuracy. More appropriate settings that balance accuracy and efficiency shall be discussed in the next section. On the collocation discretisation, we use 2500 and 400 constant elements for the part of the wall $\Gamma_0$ and the scatterer surface $\partial \Omega_0$, respectively. The layer potentials are evaluated with the $10$ points Gaussian quadrature (for both the proposed and conventional BEMs), while the Fourier integrals are evaluated by the trapezoidal rule, in which 40 integral points are used per unit wavelength of the shortest wave. The Burton-Miller parameter $\beta$ is set as either $\beta=0$, $\mathrm{i}/k$, or $-\mathrm{i}/k$. 

\begin{figure}[h]
 \centering
 \includegraphics[scale=0.55]{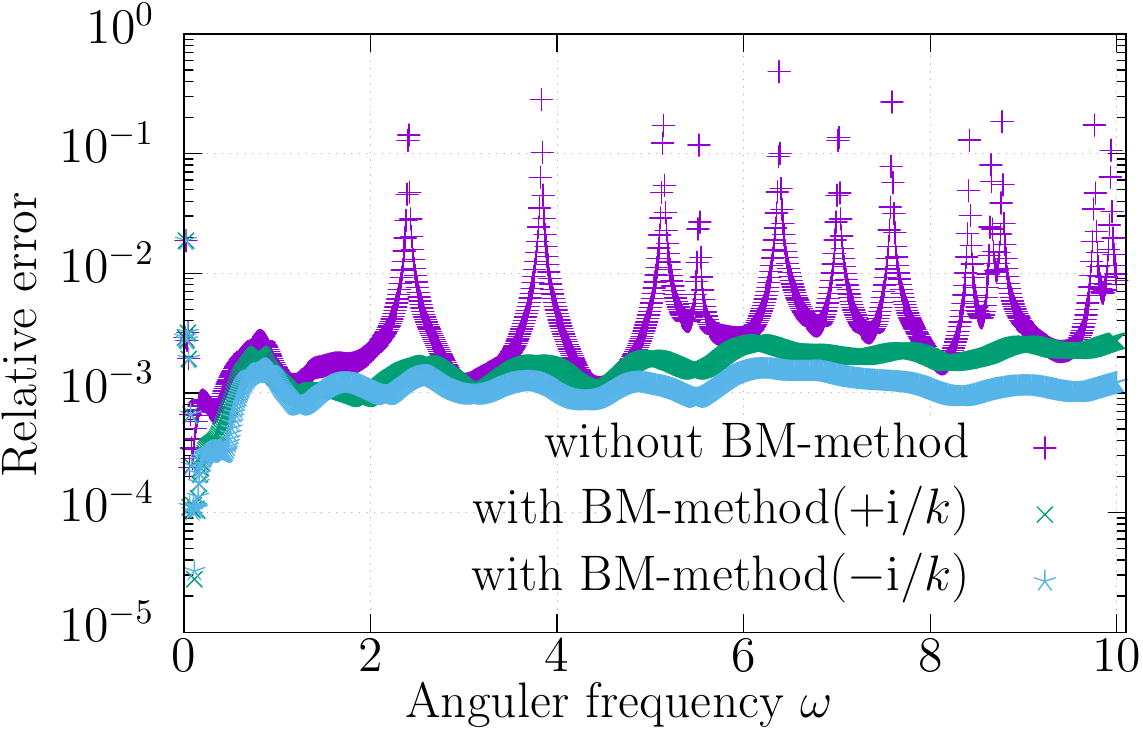}
 \caption{Relative error between the results by the proposed method and the method of image.}\label{fig:pmBM_noBMlogy}
\end{figure}

Figure \ref{fig:pmBM_noBMlogy} shows the relative error \eqref{eq:Relativeerror} versus the angular frequency. The relative error for the case of $\beta=0$ (labelled as ``without BM-method'') exhibits extremely large value at several angular frequencies. Since those frequencies correspond to the zeros of the Bessel functions, we may conclude that the poor accuracy is caused by the fictitious eigenvalue problem due to the single layer ansatz \eqref{eq:bmintegral} with $\beta=0$. On the contrary, the relative errors with Burton-Miller's method with $\beta=\pm \mathrm{i}/k$ (labelled as ``with BM-method ($\pm \mathrm{i}/k$)'') avoid the fictitious eigenvalues. In addition, the proposed Burton-Miller formulation not only avoids the fictitious eigenvalue problems but also gives smaller relative errors than the original formulation in a broad range of frequencies. It might also be concluded that the method with $\beta=-\mathrm{i}/k$ outperforms that with $\beta=\mathrm{i}/k$ from the figure. 

\subsection{Parameter sweep for $M_0$ and $N_0$} \label{sec:M0N0}
In the hybrid BEM, the choice of parameters $M_0$ and $N_0$ is of significant importance as it directly influences both the efficiency and accuracy of the method. Since these parameters are utilised to truncate the infinite integral intervals (see \eqref{eq:HBsol0} and \eqref{eq:HBsol}), larger values naturally result in more accurate solutions. Considering the efficiency aspect, on the other hand, it is desirable to keep these parameters as small as possible in order to reduce the number of required integral points for the integrals. We are thus motivated to determine the optimal values for these parameters. In this section, we solve several benchmark problems with various $M_0$ and $N_0$, aiming to find the optimal settings that provide sufficiently accurate results. As the accuracy indicator, we here use the $L_1$-relative errors in real and imaginary parts of $u|_{\partial\Omega_0}$. The reference solution is obtained by the method of image.

\begin{figure}[h]
 \centering
 \includegraphics[scale=0.48]{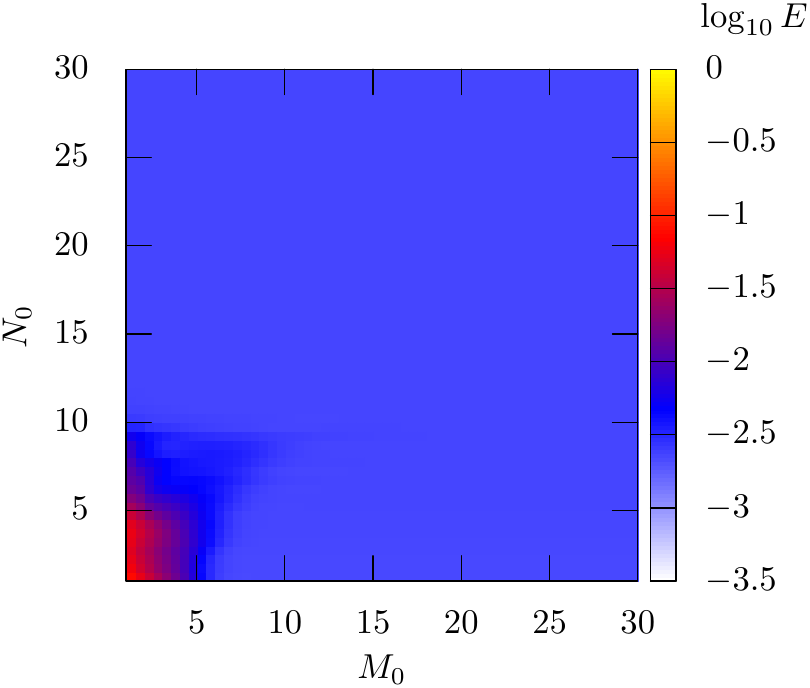}
 \includegraphics[scale=0.48]{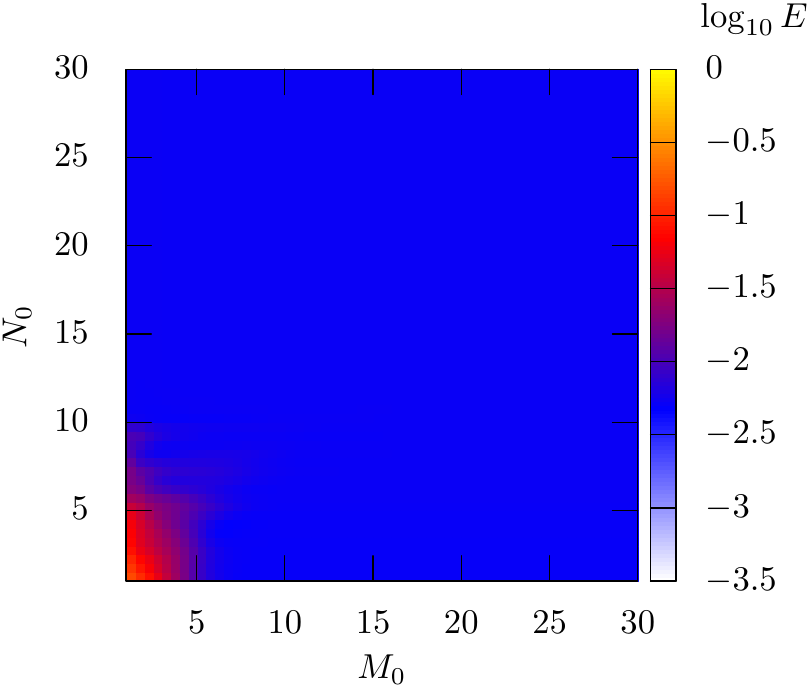}
 \caption{The relative errors $E$ in the real (left) and imaginary (right) parts of the boundary trace of $u$.}\label{fig:MNmap_10.0base}
\end{figure}

As the first example, we use the same configuration as the one in the previous example, i.e. $\Omega_0$ is the unit circle located at $(0.0, 3.0)$, and the source is located at $(1.0, 3.0)$ vibrating with $\omega=10.0$. Here, we put boundary elements of equal length on the truncated wall boundary $\Gamma_0$ and the scatterer surface $\partial \Omega_0$ and set the length of each element as the 1/40 of the wavelength for both the proposed and conventional methods. Figure \ref{fig:MNmap_10.0base} shows the relative errors in the real and imaginary parts. As expected, the combination of large $M_0$ and $N_0$ gives sufficiently accurate result. One observes declined accuracy in the case that $M_0$ and $N_0$ are smaller than 5 for this setting. 

\begin{figure}[h]
 \centering
 \includegraphics[scale=0.48]{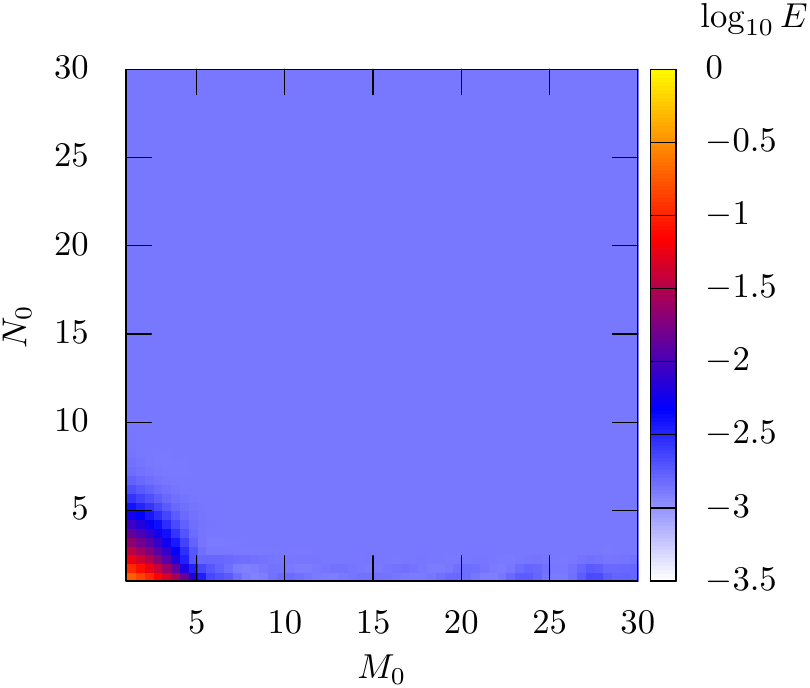}
 \includegraphics[scale=0.48]{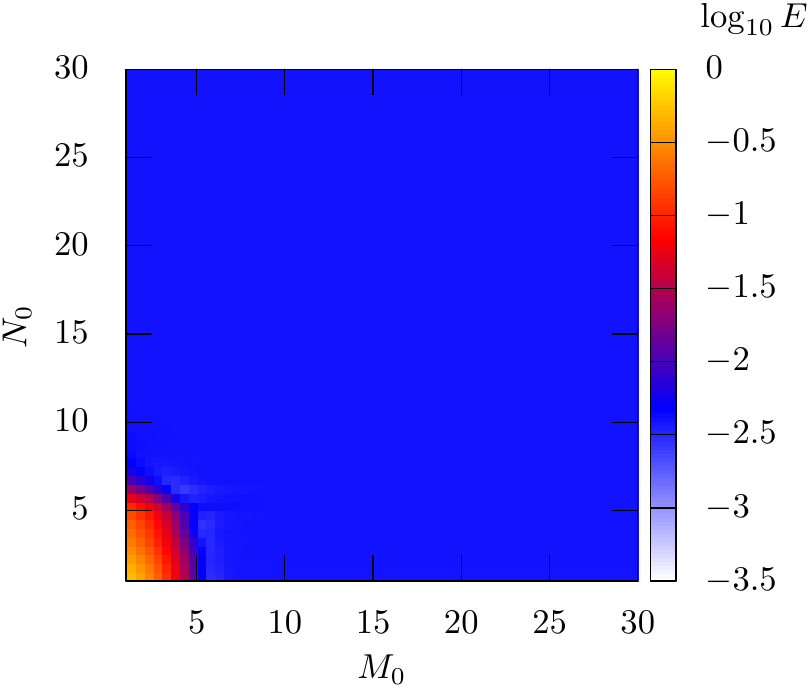}
 \includegraphics[scale=0.48]{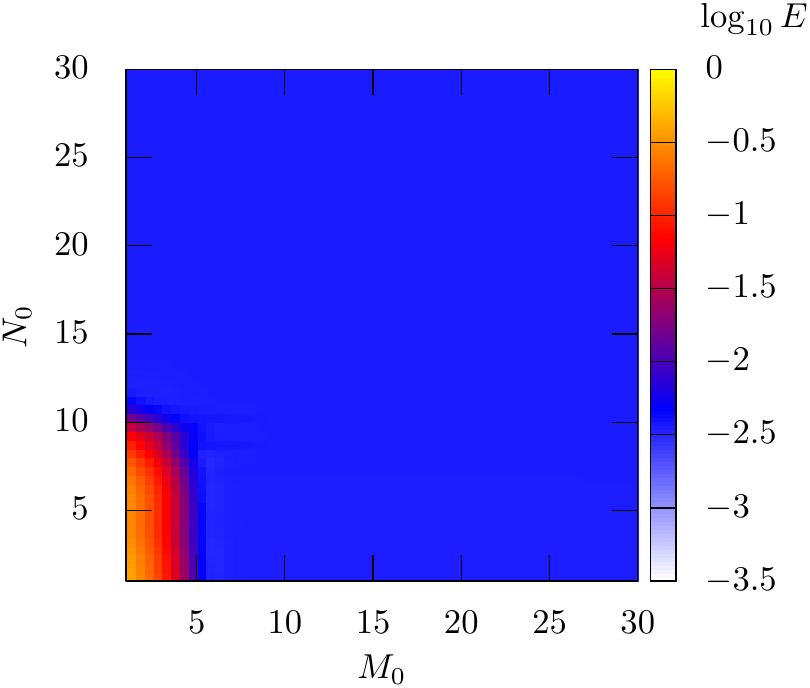}
 \caption{The relative errors $E$ in the real part of the boundary trace of $u$ in the case that the source is close to $\Gamma$. The left, centre, and right figures respectively show the results with $\omega=2.0$, $6.0$, and $10.0$.}\label{fig:MNmapreal}
\end{figure}

To assess the performance of the proposed method with more severe situation, we then position the source at $(0.0, 0.010)$ in close proximity to the infinite boundary $\Gamma$. Here, to ensure accuracy, we use a graded mesh, wherein shorter boundary elements are concentrated near the source, while the length of the longest elements is 1/40 of the wavelength for a wave oscillating with $\omega=10.0$. With this example, we also check the relation between the parameters and the underlying wavelength. We here run the computations with the angular frequency set as $\omega=2.0, 6.0,$ or $10.0$. The other configurations are all set as the previous example. Figure \ref{fig:MNmapreal} shows the results. Comparing the rightmost one with \figref{fig:MNmap_10.0base}, the red area where the accuracy is declined expanded in this case. None the less, we observe that satisfactory results are obtained if the $M_0$ and $N_0$ are set larger respectively than 5 and 10. Figure \ref{fig:MNmapreal} also indicates that the change in the angular frequency $\omega$ has a strong influence on the $N_0$ setting, which is considered to be reasonable because this parameter truncates the Fourier integral in \eqref{eq:sec2_4_2}. From these observations, we may conclude that $N_0$ should be set larger than the angular frequency $\omega$ of interest. 

\begin{figure}[h]
 \centering
 \includegraphics[scale=0.48]{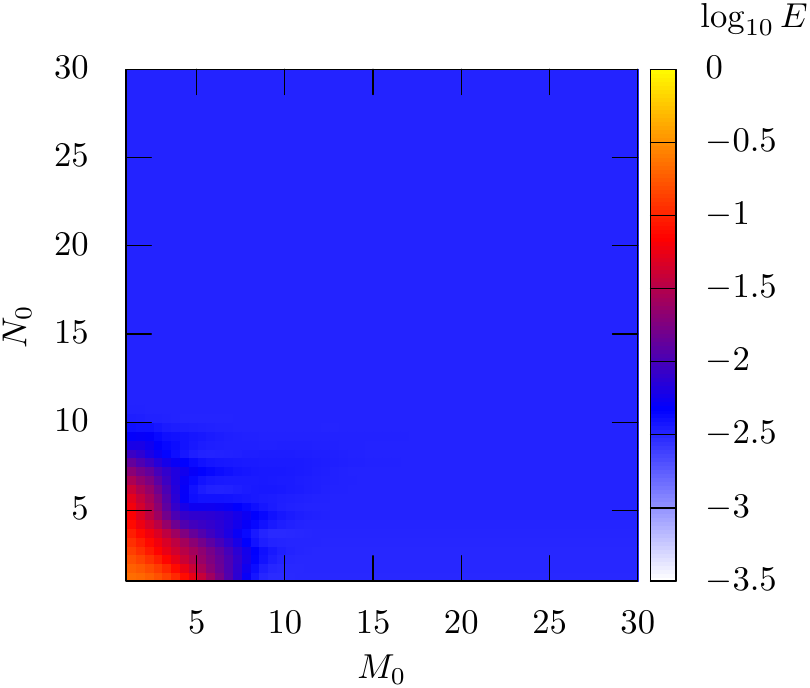}
 \includegraphics[scale=0.48]{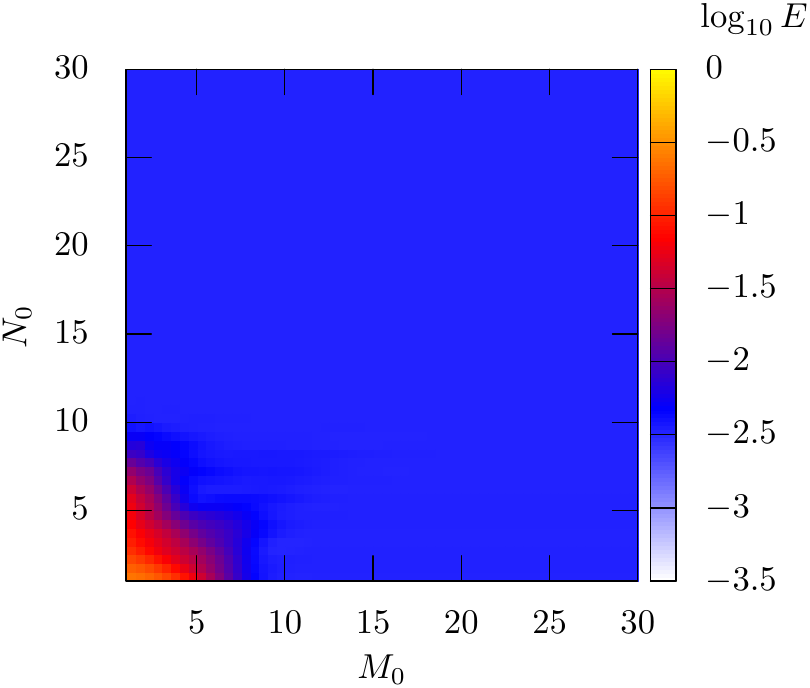}
 \caption{The relative errors $E$ in the real (left) and imaginary (right) parts of the boundary trace of $u$ in the case with a large scatterer.}\label{fig:MNmapbig}
\end{figure}

Our final example focuses on investigating the impact of scatterer size and the source location on the appropriate setting for $M_0$. In this case, we introduce a circular scatterer with the radius of 2.0 and centre at (0.0, 2.5). The source is now set at (2.0, 4.5). All other settings remain the same as in the first example. Figure \ref{fig:MNmapbig} illustrates the relative errors in the real and imaginary parts of $u$, from which we confirm that, as can easily be imagined, a larger value for $M_0$ is required to ensure the accuracy when dealing with a larger scatterer. Also, the position of the source may have a strong influence on $M_0$. From these observations, we may conclude that $M_0$ needs to be set several times the width of the rectangular enclosing the scatterer and the source. 

\section{Hybrid BEM for cavity scattering}
\subsection{Settings}
\begin{figure}[h]
 \centering
 \includegraphics[scale=0.23]{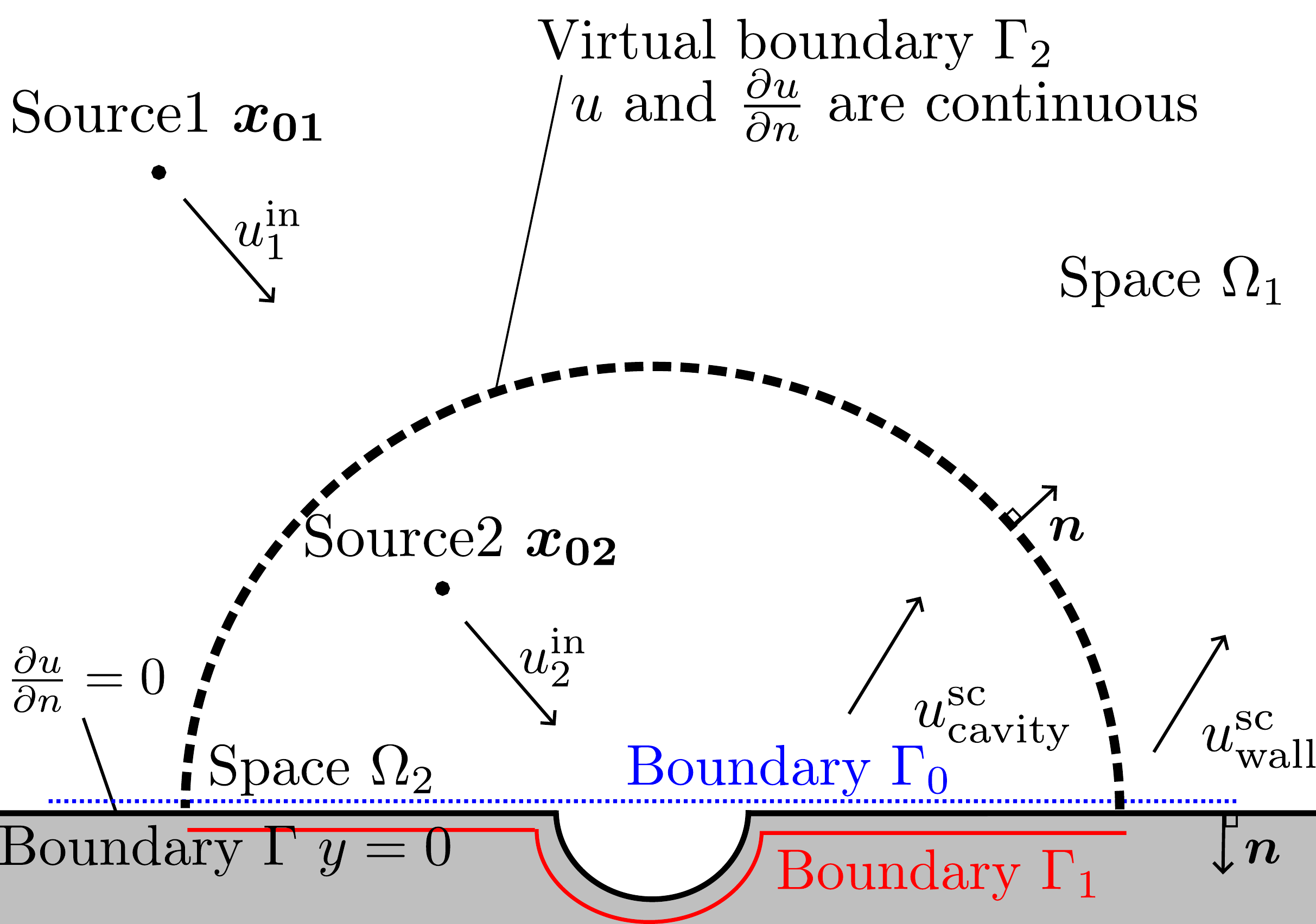}
 \caption{Cavity scattering.}\label{fig:cavity1}
\end{figure}

We here consider the acoustic scattering in a semi-infinite domain $\Omega$ with a locally perturbed boundary $\Gamma$ as depicted in \figref{fig:cavity1}.
Here, the unperturbed part of $\Gamma$ is assumed on the line $y=0$ without loss of generality. Since the perturbed domain (called cavity in this paper) is concaved into the region $y<0$, the scattering problem cannot be solved by the conventional method of image. To address the problem, we here introduce a bounded subdomain $\Omega_2$ of $\Omega$ that includes the cavity. The outside of $\Omega_2$ is denoted as $\Omega_1:=\Omega\setminus\overline{\Omega_2}$. The boundary of $\Omega_2$ is partitioned into $\Gamma_1:=\partial\Omega_2\cap \Gamma$ and $\Gamma_2:=\partial\Omega_2\setminus\overline{\Gamma_1}$, and $\Gamma_2$ is henceforth called the virtual boundary in this paper. The normal vector on $\Gamma$ is naturally defined as positive when it is directed from $\Omega$. We define, for convenience, the unit normal also on the virtual boundary $\Gamma_2$, which is directed from $\Omega_2$. The point sound sources may be allocated both in $\Omega_1$ and $\Omega_2$. The source position in $\Omega_i$ is denoted as $\bs{x_{0i}}$. Although a single source is assumed in each subdomain in the figure and the formulations to follow, it is straightforward to modify the formulation to delete one of those and/or add some more sources. On the physical boundary $\Gamma$, the homogeneous Neumann boundary condition is again imposed. On $\Gamma_2$, since this is virtual, the sound pressure and its normal gradient should be continuous as
\begin{align}
 \lim_{\epsilon \,\downarrow\, 0} u(\bs{x}+\epsilon\bs{n}(\bs{x}))&=\lim_{\epsilon\,\downarrow\, 0} u(\bs{x}-\epsilon\bs{n}(\bs{x})) ~~~~\bs{x} ~\mathrm{on} ~\Gamma_2, \label{4_consis1}\\
 \lim_{\epsilon \,\downarrow\, 0} \frac{\partial u(\bs{x}+\epsilon\bs{n}(\bs{x}))}{\partial \bs{n}(\bs{x})}&=\lim_{\epsilon \,\downarrow\, 0} \frac{\partial u(\bs{x}-\epsilon\bs{n}(\bs{x})) }{\partial \bs{n}(\bs{x})}~~\bs{x} ~\mathrm{on} ~\Gamma_2. \label{4_consis2}
\end{align}

\subsection{Formulation}
To solve the cavity scattering problem defined in the preceding section, we assume that the total field has the following representation respectively in each subdomain:
\begin{align}
u(\bs{x})&=u^{\mathrm{in}}_1(\bs{x})+u^{\mathrm{sc}}_\mathrm{wall}(\bs{x})+u_{\mathrm{virtual}}(\bs{x}) ~~~~\bs{x} \in \Omega_1, \label{4_u-1}\\  u(\bs{x})&=u^{\mathrm{in}}_2(\bs{x})+u^{\mathrm{sc}}_{\mathrm{cavity}}(\bs{x})+u_{\mathrm{virtual}}(\bs{x})~~\bs{x} \in \Omega_2, \label{4_u-2}
\end{align}
where $u^\mathrm{in}_i(\bs{x}):=G_k(\bs{x}, \bs{x_{0i}})$ is the incident field induced by the source $\bs{x_{0i}}\in\Omega_i$, and the integral representation for the scattered fields from the wall, cavity and virtual boundary are respectively defined as
\begin{align}
 u^{\mathrm{sc}}_{\mathrm{wall}}(\bs{x})&=S_{\Gamma_0}[\sigma_W](\bs{x}) +F_{I_0}[ \widehat{\xi_W}](\bs{x}), \label{18} \\
 u^{\mathrm{sc}}_{\mathrm{cavity}}(\bs{x})&=S_{ \Gamma_1}[\sigma](\bs{x}),\\
 u_{\mathrm{virtual}}(\bs{x})&=S_{ \Gamma_2}[\sigma](\bs{x})+ D_{\Gamma_2}[\mu](\bs{x}),
\end{align}
with the unknown density functions $\sigma$, $\mu$, and $\widehat{\xi_W}$. It should be noted that $\Gamma_0=\{\bs{x}\mid -M_0<x<M_0,~y=0\}$ in \eqref{18} may not necessarily be the subset of $\Gamma$ in this case (see \figref{fig:cavity1}). Some parts of $\Gamma_0$ can thus also be virtual. For the definitions of the single layer potential $S_{\Gamma_i}~(i=0,~1,~2)$, the double layer potential $D_{\Gamma_2}$ and the Sommerfeld integral $F_{I_0}$, see \eqref{eq:sec2_4_1}, \eqref{7}, and \eqref{eq:sec2_4_2}, respectively. Also, $\sigma_W$ is the density function $\sigma$ multiplied by the window function; see \eqref{eq:sgmW} and \eqref{eq:window}. 

We then derive the boundary integral equations for determining the density functions. Corresponding to \eqref{eq:H1BM_1}, we set the following integral equation for $\bs{x}\in\Gamma_0$: 
\begin{align}
  \frac{1}{2} \sigma|_{\Gamma_0}(\bs{x}) +D^*_{\Gamma_0}[\sigma](\bs{x})+D^*_{\Gamma_2}[\sigma](\bs{x})+N_{\Gamma_2}[\mu](\bs{x})=-\frac{\partial G_k(\bs{x},\bs{x_{01}})}{\partial \bs{n}(\bs{x})}. \label{eq:4_bc1}
\end{align}
Since the boundary integral equation \eqref{eq:4_bc1} is stemmed from the exterior integral representation \eqref{4_u-1} for $\bs{x}\in\Omega_1$, only the layer potentials from $\Gamma_0$ and the virtual boundary $\Gamma_2$, as well as the point source located in the outer region $\Omega_1$ are involved. Since the integral equation \eqref{eq:4_bc1} does not take into account the scattering from $\Gamma\setminus\overline{\Gamma_0}$, we need additionally the following equation in the spectral space:
\begin{equation}
  \begin{split}
      \frac{1}{2\pi} \int^{\infty}_{-\infty} \frac{1}{2}\widehat{\sigma_W} e^{i \lambda x}d\lambda +\frac{1}{2\pi} \int^{\infty}_{-\infty}  \widehat{D^*_{\Gamma_0}[\sigma]}  e^{i \lambda x}d\lambda +\frac{1}{2\pi} \int^{\infty}_{-\infty}  \widehat{D^*_{\Gamma_2}[\sigma]}  e^{i \lambda x}d\lambda \\
      +\frac{1}{2\pi} \int^{\infty}_{-\infty}  \widehat{N_{\Gamma_2}[\mu]}e^{i \lambda x}d\lambda +H_{I_0}[\widehat{\xi_W}]
      =-\frac{1}{2\pi} \int^{\infty}_{-\infty} \widehat{ \frac{\partial G_k}{\partial \bs{n}(\bs{x})}} e^{i \lambda x}d\lambda,
  \end{split}
\label{22}
\end{equation}
in order to make the total field satisfy the boundary condition on the whole $\Gamma$. In \eqref{22}, the finite integral interval $I_0=(-N_0,~N_0)$ for $H_{I_0}$ can be replaced by $(-\infty,~\infty)$ due to the local supportness of the integrand (see \eqref{eq:H_I_0}). As a result, we have
\begin{equation}
 \begin{split}
  \frac{1}{2}\widehat{\sigma_W}  + \widehat{D^*_{\Gamma_0}[\sigma]}   + \widehat{D^*_{\Gamma_2}[\sigma]}
  + \widehat{N_{\Gamma_2}[\mu]} +\frac{1}{2} \widehat{\xi_W}
  =-\widehat{ \frac{\partial G_k}{\partial \bs{n}(\bs{x})}}. \label{eq:4_bc2}
 \end{split}
\end{equation}
We then consider the boundary condition on the cavity surface $\Gamma_1$. Taking the limit of (the normal derivative of) the integral representation \eqref{4_u-2} as $\Omega_1\ni\bs{x}\rightarrow \bs{x}\in \Gamma_1$, we have
 \begin{align}
   \frac{1}{2} \sigma|_{\Gamma_1}(\bs{x}) +D^*_{\Gamma_1}[\sigma](\bs{x})+D^*_{\Gamma_2}[\sigma](\bs{x})+N_{\Gamma_2}[\mu](\bs{x})=-\frac{\partial G_k(\bs{x},\bs{x_{02}})}{\partial \bs{n}(\bs{x})}, \label{eq:4_bc3}
\end{align}
for $\bs{x}\in\Gamma_2$. We then consider the boundary conditions on the virtual boundary $\Gamma_2$. By taking the appropriate limits for the integral representations \eqref{4_u-1} and \eqref{4_u-2} and substituting the resulting traces into the boundary conditions \eqref{4_consis1} and \eqref{4_consis2}, we have the following boundary integral equations: 
\begin{align}
  S_{\Gamma_0}[\sigma](\bs{x}) +F_{I_0}[\widehat{\xi_W}](\bs{x}) +\mu(\bs{x}) -S_{\Gamma_1}[\sigma](\bs{x})&= G_k(\bs{x},\bs{x_{02}})- G_k(\bs{x},\bs{x_{01}}), \label{eq:4_bc4}\\
  D^{*}_{\Gamma_1}[\sigma](\bs{x})+\sigma|_{\Gamma_2}(\bs{x}) -D^{*}_{\Gamma_0}[\sigma](\bs{x}) -H_{I_0}[\widehat{\xi_W}](\bs{x})&=\frac{\partial G_k(\bs{x},\bs{x_{01}})}{\partial \bs{n}(\bs{x})}-\frac{\partial G_k(\bs{x},\bs{x_{02}})}{\partial \bs{n}(\bs{x})}, \label{eq:4_bc5}
\end{align}
for $\bs{x}\in\Gamma_2$. The equations thus obtained \eqref{eq:4_bc1}, \eqref{eq:4_bc2}, \eqref{eq:4_bc3}, \eqref{eq:4_bc4}, and \eqref{eq:4_bc5} are simultaneously solved after discretisation to determine the unknown density functions.

\begin{figure}[h]
 \centering
 \includegraphics[scale=0.3]{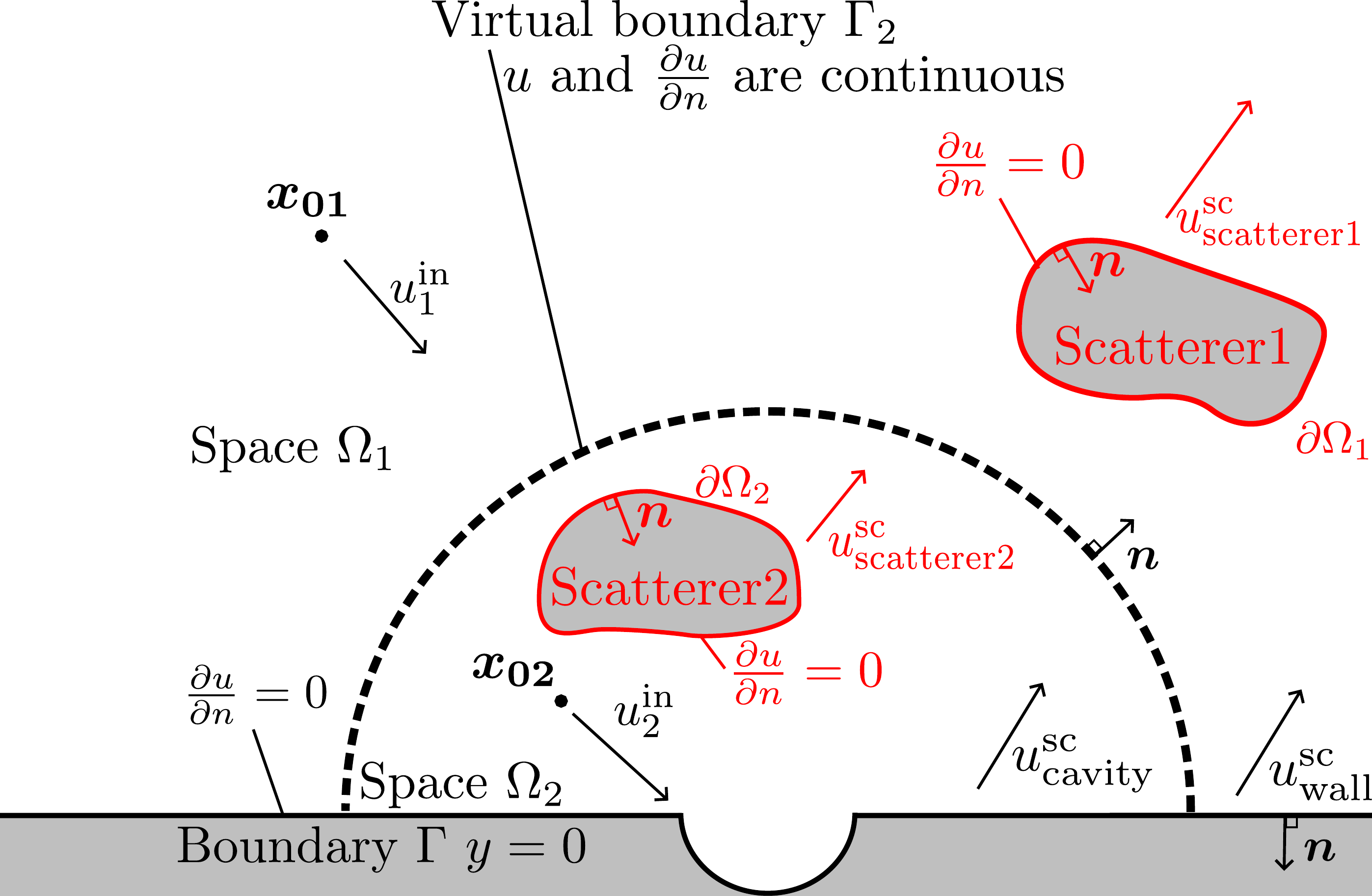}
 \caption{The setting of a cavity scattering problem with some scatterers. }\label{fig:cavity2}
\end{figure}

In a similar way, we may solve the scattering problems when several rigid scatterers exist in a half-space with a cavity. We here show in \figref{fig:cavity2} an illustrative sketch of the situation, in which a scatterer exists in each subdomain $\Omega_1$ and $\Omega_2$, each of which is labelled as Scatterer1 and Scatterer2 in the figure. 
In this case, we need to modify the integral representations \eqref{4_u-1} and \eqref{4_u-2} respectively as
\begin{align}
  u(\bs{x})&=u^{\mathrm{in}}_1(\bs{x})+u^{\mathrm{sc}}_\mathrm{wall}(\bs{x})+u^{\mathrm{sc}}_{\mathrm{scatterer1}}+u_{\mathrm{virtual}}(\bs{x}) ~~~~\bs{x} \in \Omega_1, \label{5_u-1}\\
  u(\bs{x})&=u^{\mathrm{in}}_2(\bs{x})+u^{\mathrm{sc}}_{\mathrm{cavity}}(\bs{x})+u^{\mathrm{sc}}_{\mathrm{scatterer2}}+u_{\mathrm{virtual}}(\bs{x})~~\bs{x} \in \Omega_2, \label{5_u-2}
\end{align}
in which the scattered field from the scatterer in $\Omega_1$ and $\Omega_2$ is respectively denoted as $u^\mathrm{sc}_\mathrm{scatterer1}$ and $u^\mathrm{sc}_\mathrm{scatterer2}$. Considering the fictitious eigenvalue problem, their integral representations are respectively given as
\begin{align}
  u^{\mathrm{sc}}_{\mathrm{scatterer1}}(\bs{x})=S_{ \partial \Omega_1}[\sigma](\bs{x})+\beta D_{ \partial \Omega_1}[\sigma](\bs{x}),\\
  u^{\mathrm{sc}}_{\mathrm{scatterer2}}(\bs{x})=S_{ \partial \Omega_2}[\sigma](\bs{x})+\beta D_{ \partial \Omega_2}[\sigma](\bs{x}), 
\end{align}
where $\sigma$ is the unknown density function, and $\beta\in\mathbb{C}$ is the parameter for the Burton-Miller method. The unit normal $\bs{n}$ on $\partial \Omega_1$ and $\partial \Omega_2$ is directed from $\Omega$. The boundary integral equations for the unknown density functions can be obtained as the previous case. Here, we just list the results as follows:
\begin{equation}
  \begin{split}
      \frac{1}{2} \sigma|_{\Gamma_0}(\bs{x}) +D^*_{\Gamma_0}[\sigma](\bs{x})+D^*_{\partial \Omega_1}[\sigma](\bs{x})+ \beta N_{\partial \Omega_1}[\sigma](\bs{x})+D^*_{\Gamma_2}[\sigma](\bs{x})+N_{\Gamma_2}[\mu](\bs{x})\\
      =-\frac{\partial G_k(\bs{x},\bs{x_{01}})}{\partial \bs{n}(\bs{x})},\label{eq:5_bc1}
  \end{split}
\end{equation}
\begin{equation}
  \begin{split}
       \frac{1}{2}\widehat{\sigma_W}  + \widehat{D^*_{\Gamma_0}  [\sigma]}
       + \widehat{D^*_{\partial \Omega_1}[\sigma]}
       +\beta \widehat{N_{\partial \Omega_1}[\sigma]}
       + \widehat{D^*_{\Gamma_2}[\sigma]}
       + \widehat{N_{\Gamma_2}[\mu]} 
       +\frac{1}{2} \widehat{\xi_W}
      =-\widehat{ \frac{\partial G_k}{\partial \bs{n}(\bs{x})}},  \label{eq:5_bc2}
  \end{split}
\end{equation}
\begin{equation}
  \begin{split}
      D^*_{\Gamma_0}[\sigma](\bs{x})+H_{I_0}[\widehat{\xi_W}](\bs{x})+\frac{1}{2} \sigma|_{\partial \Omega_1}(\bs{x}) +D^*_{\partial \Omega_1}[\sigma](\bs{x})+ \beta N_{\partial \Omega_1}[\sigma](\bs{x})\\
      +D^*_{\Gamma_2}[\sigma](\bs{x})+N_{\Gamma_2}[\mu](\bs{x})
      =-\frac{\partial G_k(\bs{x},\bs{x_{01}})}{\partial \bs{n}(\bs{x})},\label{eq:5_bc6}
  \end{split}
\end{equation}
\begin{equation}
  \begin{split}
      \frac{1}{2} \sigma|_{\Gamma_1}(\bs{x}) +D^*_{\Gamma_1}[\sigma](\bs{x})+D^*_{\partial \Omega_2}[\sigma](\bs{x})+\beta N_{\partial \Omega_2}[\sigma](\bs{x})+D^*_{\Gamma_2}[\sigma](\bs{x})+N_{\Gamma_2}[\mu](\bs{x})\\
      =-\frac{\partial G_k(\bs{x},\bs{x_{02}})}{\partial \bs{n}(\bs{x})},\label{eq:5_bc3}
  \end{split}
\end{equation}
\begin{equation}
  \begin{split}
      D^*_{\Gamma_1}[\sigma](\bs{x})+\frac{1}{2} \sigma|_{\partial \Omega_2}(\bs{x}) +D^*_{\partial \Omega_2}[\sigma](\bs{x})+\beta N_{\partial \Omega_2}[\sigma](\bs{x})+D^*_{\Gamma_2}[\sigma](\bs{x})+N_{\Gamma_2}[\mu](\bs{x})\\
      =-\frac{\partial G_k(\bs{x},\bs{x_{02}})}{\partial \bs{n}(\bs{x})},\label{eq:5_bc7}
  \end{split}
\end{equation}
\begin{equation}
  \begin{split}
      S_{\Gamma_0}[\sigma](\bs{x}) +F_{I_0}[\widehat{\xi_W}](\bs{x})+S_{\partial \Omega_1}[\sigma](\bs{x})+\beta D_{\partial \Omega_1}[\sigma](\bs{x})
      -S_{\partial \Omega_2}[\sigma](\bs{x})-\beta D_{\partial \Omega_2}[\sigma](\bs{x}) \\
      +\mu(\bs{x}) -S_{\Gamma_1}[\sigma](\bs{x})= G_k(\bs{x},\bs{x_{02}})- G_k(\bs{x},\bs{x_{01}}), \label{eq:5_bc4}
  \end{split}
\end{equation}
\begin{equation}
  \begin{split}
  D^{*}_{\Gamma_1}[\sigma](\bs{x})-D^{*}_{\partial \Omega_1}[\sigma](\bs{x}) -\beta N_{\partial \Omega_1}[\sigma](\bs{x})+D^{*}_{\partial \Omega_2}[\sigma](\bs{x}) +\beta N_{\partial \Omega_2}[\sigma](\bs{x})\\
  +\sigma|_{\Gamma_2}(\bs{x}) -D^{*}_{\Gamma_0}[\sigma](\bs{x}) -H_{I_0}[\widehat{\xi_W}](\bs{x})=\frac{\partial G_k(\bs{x},\bs{x_{01}})}{\partial \bs{n}(\bs{x})}-\frac{\partial G_k(\bs{x},\bs{x_{02}})}{\partial \bs{n}(\bs{x})}.\label{eq:5_bc5}
\end{split}
\end{equation}

\subsection{Validation and numerical examples for the cavity scattering}
In this subsection, we first validate our hybrid BEM for cavity scattering. We then demonstrate some numerical examples of scattering by small cavities that may work as open resonators \cite{kao2008maximization}.

We may use the existing methods~\cite{lai2014fast,bao2016stability} to check the accuracy of the present method. In this paper, however, we use a simpler way for the validation; we check how accurately our results satisfy the boundary conditions. To this end, we compute the following numerical derivative of the total field: 
\begin{align}
  \frac{\partial u (\bs{x})}{\partial \bs{n}(\bs{x})} \approx \frac{u(\bs{x})-u(\bs{x}-\epsilon \bs{n}(\bs{x}))}{\epsilon}, \label{eq:sabun}
\end{align}
where $\epsilon>0$ is a parameter for the numerical difference, and check if the quantity is sufficiently small. Note that the parameter $\epsilon$ should carefully be set. The smaller $\epsilon$ is, the more accurately the right-hand side of \eqref{eq:sabun} approximates the normal derivative. On the other hand,  if $\epsilon$ is too small, the numerical difference may suffer from the loss of significant digits. Since, according to the validation results in Section \ref{sec:seido1}, our numerical results may involve the relative error of order $10^{-3}$ when 40 constant elements are adopted per unit wavelength, we here set as $\epsilon=1.0 \times 10^{-3}$.

In this computation, we set the domain with a cavity as $\Omega=\mathbb{R}^2_+\cup \Omega_\mathrm{cavity}$ with $\Omega_\mathrm{cavity} = \{\bs{x} \mid x^2+y^2 < 1 \cap y<0 \}$ being a cavity. We set the virtual boundary $\Gamma_2$ as a upper half circle centred at $(0.0,0.0)$ with radius $3.0$. $\Gamma_1$ is accordingly set as $\Gamma:=\partial\Omega\setminus \{\bs{x} \mid |x|\ge 3.0 \}$. We put a single point source {oscillating with the angular frequency $\omega=10.0$} at $\bs{x_{02}}=(0.0,2.0)$ inside the virtual boundary. 

\begin{figure}[h]
 \centering
  \includegraphics[scale=0.45]{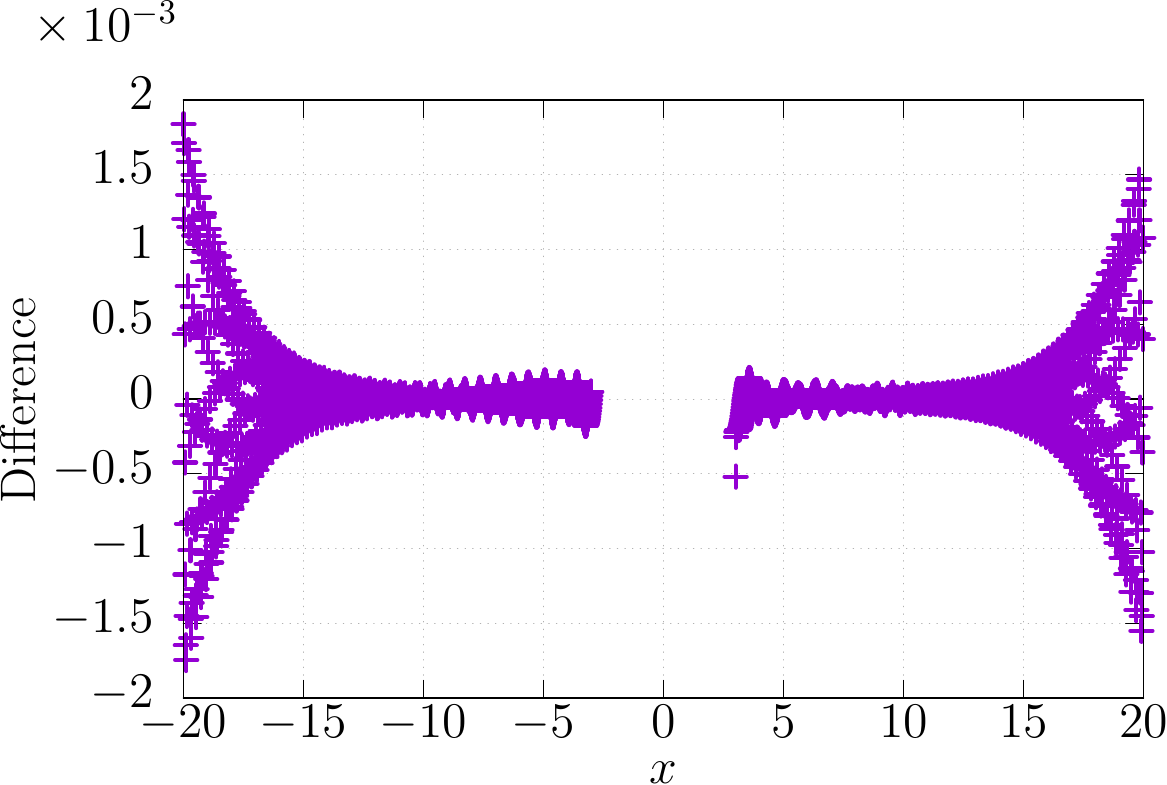}
 \caption{Numerical difference $u(\bs{x})-u(\bs{x}-\epsilon\bs{n}(\bs{x}))$ on the wall $\Gamma_0=\{\bs{x}\mid -M_0 \le x \le M_0,~y=0\}$ with $a=2$.}\label{fig:5_kekka1-1}
\end{figure}

\begin{figure}[h]
 \centering
 \includegraphics[scale=0.45]{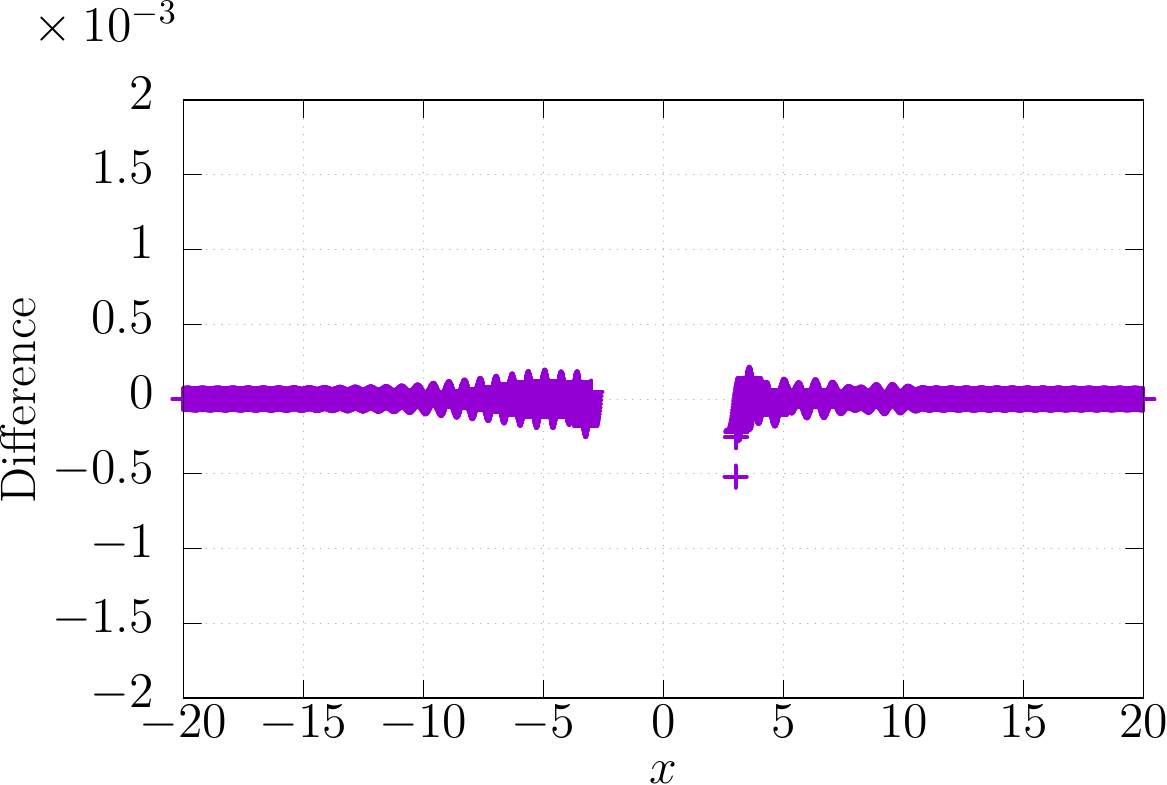}
 \caption{Numerical difference $u(\bs{x})-u(\bs{x}-\epsilon\bs{n}(\bs{x}))$ on the wall $\Gamma_0=\{\bs{x}\mid -M_0 \le x \le M_0,~y=0\}$ with $a=8$.}\label{fig:5_kekka1-3}
\end{figure}

\begin{figure}[h]
 \centering
 \includegraphics[scale=0.45]{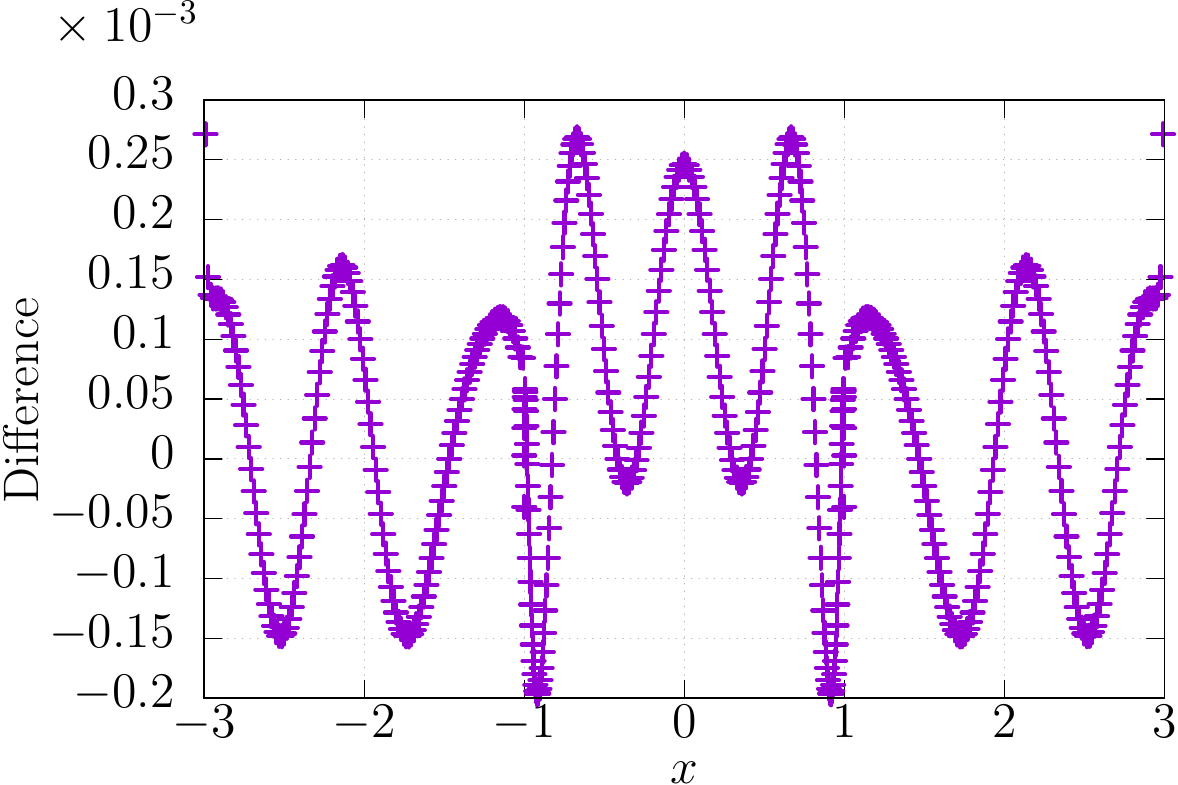}
 \caption{Value of numerical derivative on the cavity $\Gamma_1$($a=8$).}\label{fig:5_kekka1-4}
\end{figure}

\begin{figure}[h]
  \centering
  \includegraphics[scale=0.18]{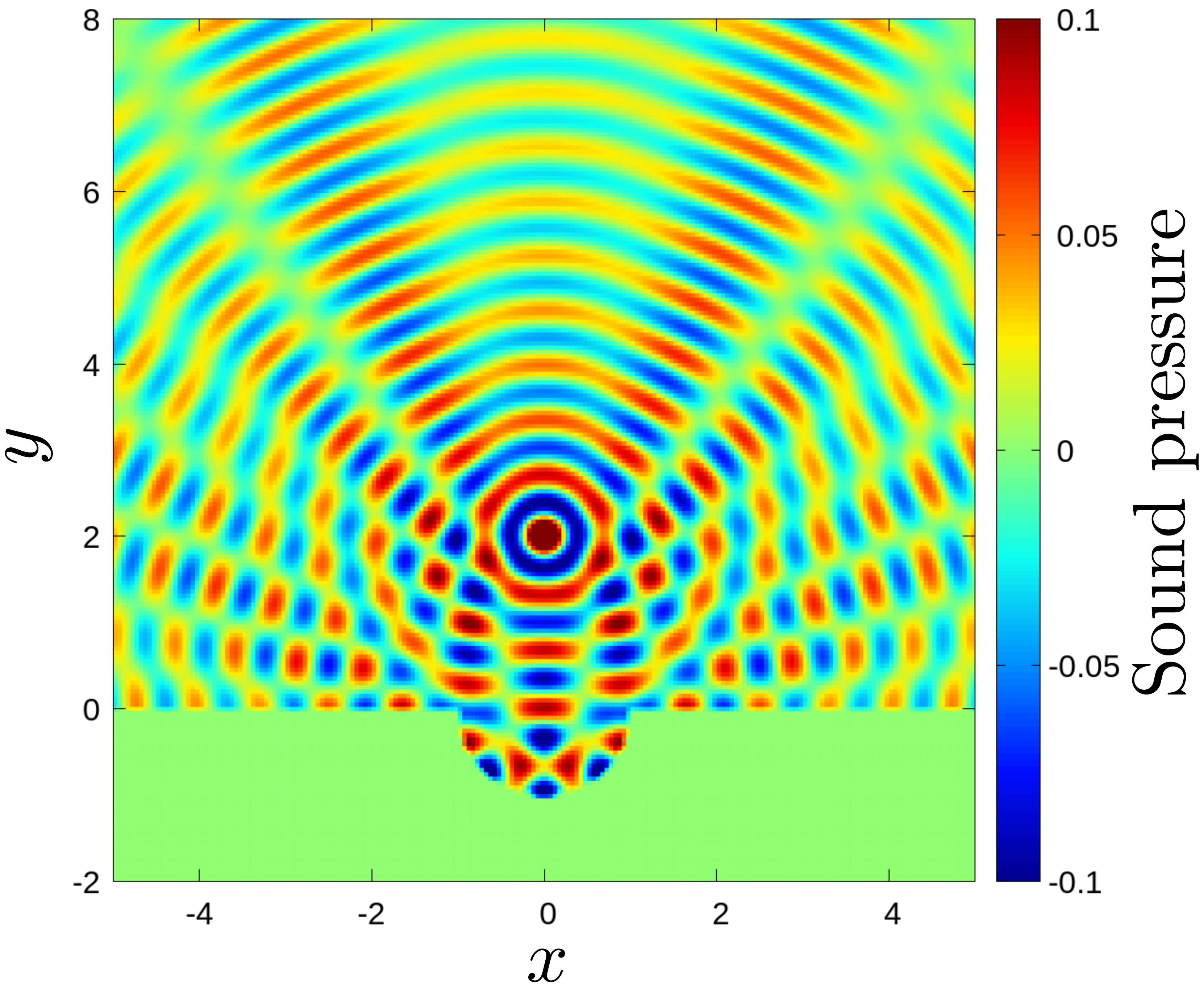}
  \caption{Real part of the total field.}\label{fig:5_kekka1-5}
\end{figure}

We use the parameters $M_0=N_0=20$ for the hybrid integrals, which is considered sufficiently large as indicated in Section \ref{sec:M0N0}. All the boundaries are discretised into piecewise constant boundary elements of length $\lambda/40$, where $\lambda$ is the underlying wavelength. The Fourier integral is evaluated by the trapezoidal rule with 40 integral points per unit wavelength of the shortest relevant wave. Figure \ref{fig:5_kekka1-1} shows the difference between $u(\bs{x})$ and $u(\bs{x}-\epsilon\bs{n}(\bs{x}))$ for $\bs{x}\in\Gamma_0$ with the integral path deformation parameter being $a=2.0$. 
One may observe that, for $|x|<15$, the difference is less than or comparable to $\epsilon$, from which one finds that the Neumann boundary condition is accurately satisfied. The difference seems, however, to diverge as $|x|$ becomes large. The poor accuracy for large $|x|$ may caused by the parameter setting for the contour deformation. This can be seen from the definition of the Sommerfeld integral \eqref{eq:Sommer} and \eqref{eq:contour_deformation}. If $a$ is small for the contour deformation \eqref{eq:contour_deformation}, the factor $e^{\mathrm{i}\lambda x}$ may diverge for the case of large $|x|$. We are thus motivated to compute again the numerical difference with larger setting for $a$, i.e. $a=8.0$. The result is summarised in \figref{fig:5_kekka1-3}.
In \figref{fig:5_kekka1-3}, the difference divergence is suppressed even for large $|x|$ and the difference is less than $\epsilon$, from which we conclude that, if the parameter is appropriately set, the proposed BEM gives numerical solution that satisfies the boundary condition on $\Gamma_0$. We then show the numerical derivative plot in the vicinity of cavity surface $\Gamma_1$ (with the parameter $a=8.0$) in \figref{fig:5_kekka1-4}, from which we confirm that the Neumann condition is satisfied on this boundary as well. 
Figure \ref{fig:5_kekka1-5} shows the real part of the total field around the cavity. 
We observe no jumps in the total field on the virtual boundary $\Gamma_2$. With these observations, we may conclude that our hybrid BEM is validated for the cavity scattering. 

\begin{figure}[h]
  \centering
  \includegraphics[scale=0.3]{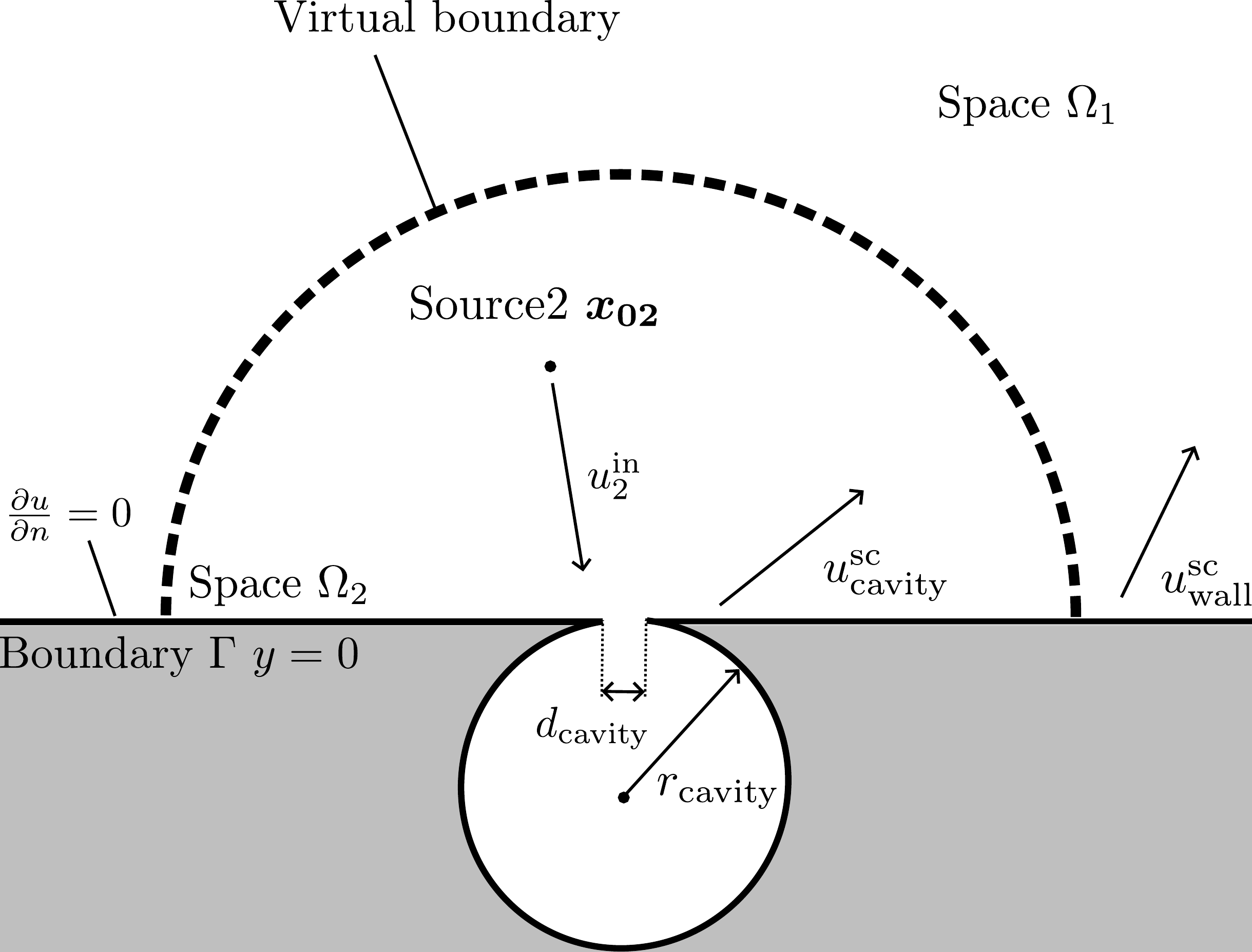}
  \caption{Acoustic scattering by the Helmholtz-resonator-like cavity}\label{fig:5_setting2}
\end{figure}

\begin{figure}[h]
  \centering
  \includegraphics[scale=0.55]{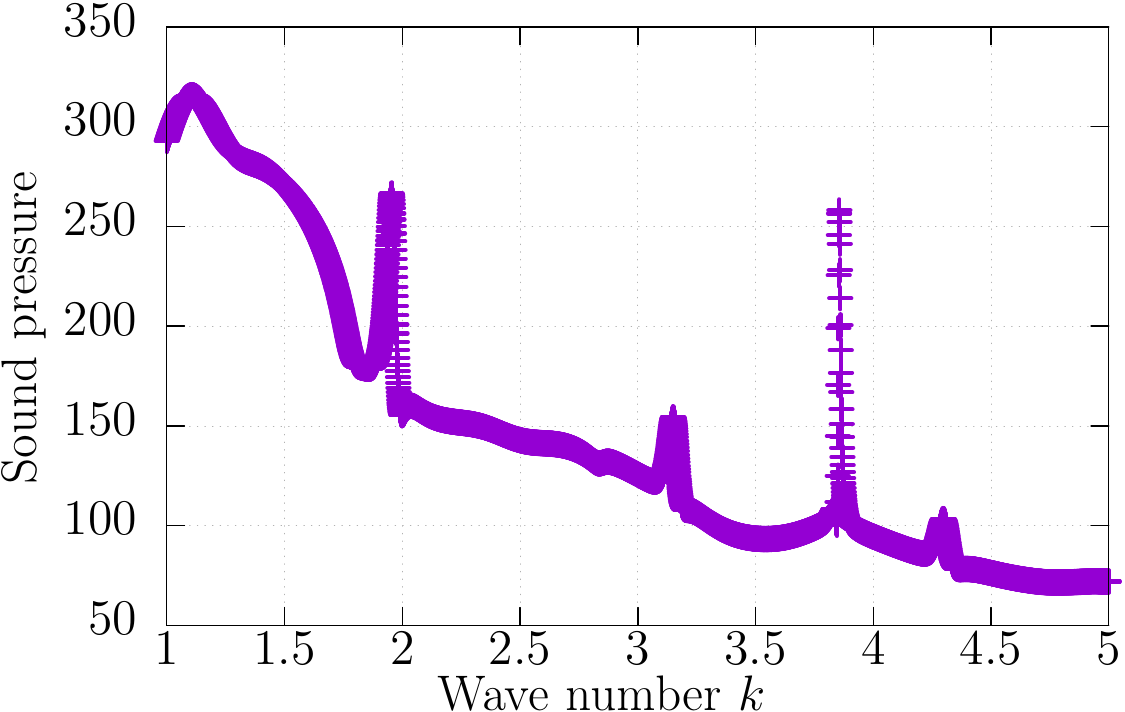}
  \caption{Frequency sweep of the sound intensity for the Helmholtz-resonator-like cavity}\label{fig:5_suml2-3}
\end{figure}

\begin{figure}[h]
  \centering
  \includegraphics[scale=0.13]{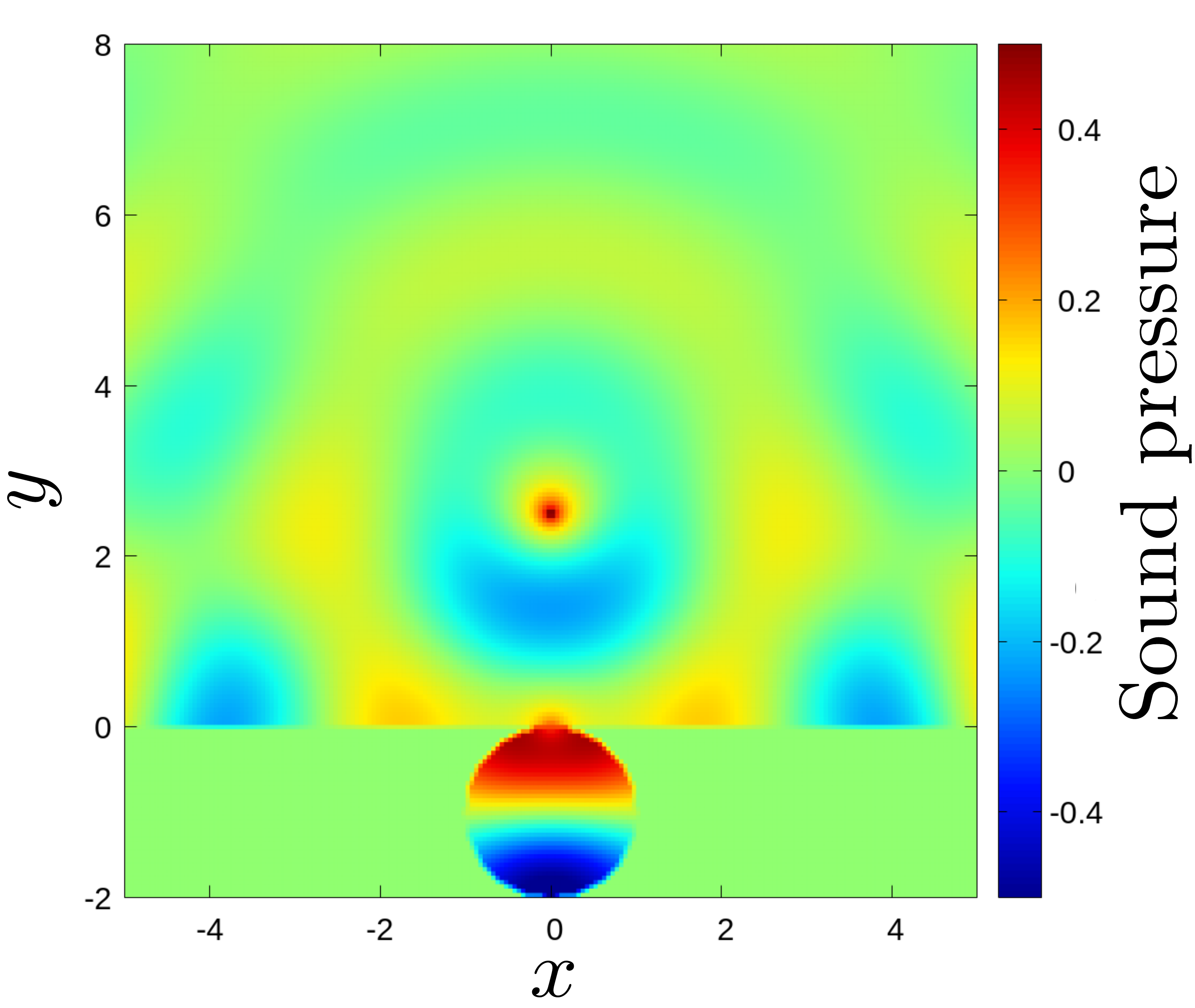}
  \caption{Sound pressure distribution with the incident wave number $k=1.960$.}\label{fig:5_reson}
\end{figure}

In what follows, we show two numerical examples related to cavity resonators. The first one is inspired by the Helmholtz resonator, and its setting is illustrated in \figref{fig:5_setting2}. 
In this case, we hollowed a circular cavity centred at $(0.0, -1.0)$ with opening $d_\mathrm{cavity}=0.1$. We set at $\bs{x_{02}}=(0.0,2.5)$ a single point source inside the virtual boundary. The rest of the settings e.g. the virtual boundary, $M_0$, $N_0$ and $a$, discretisation criteria etc are the same as the ones used in the previous example. With these settings, we sweep the frequency response of the cavity by changing the wavenumber $k$ from $1.0$ to $5.0$ by $0.001$. Figure \ref{fig:5_suml2-3} shows the result, whose vertical axis indicates the sum of the sound intensity $|u|^2$ evaluated at 40401 lattice points distributed in a rectangular domain defined as $-5.0 \le x \le 5.0$ and {$-2.1 \le y \le 8.0$}. 
The plot in \figref{fig:5_suml2-3} exhibits sharp peaks, which may correspond to the resonance in the cavity. To check the acoustic behaviour at these excitation frequencies, we show the real part of the sound pressure at $k=1.960$ in the vicinity of the cavity in \figref{fig:5_reson}.
Figure \ref{fig:5_reson} shows that the acoustic wave is localised in the cavity and thus confirms that the sharp peak in the sound intensity is caused by the resonance. 

\begin{figure}[h]
 \centering
 \includegraphics[scale=0.30]{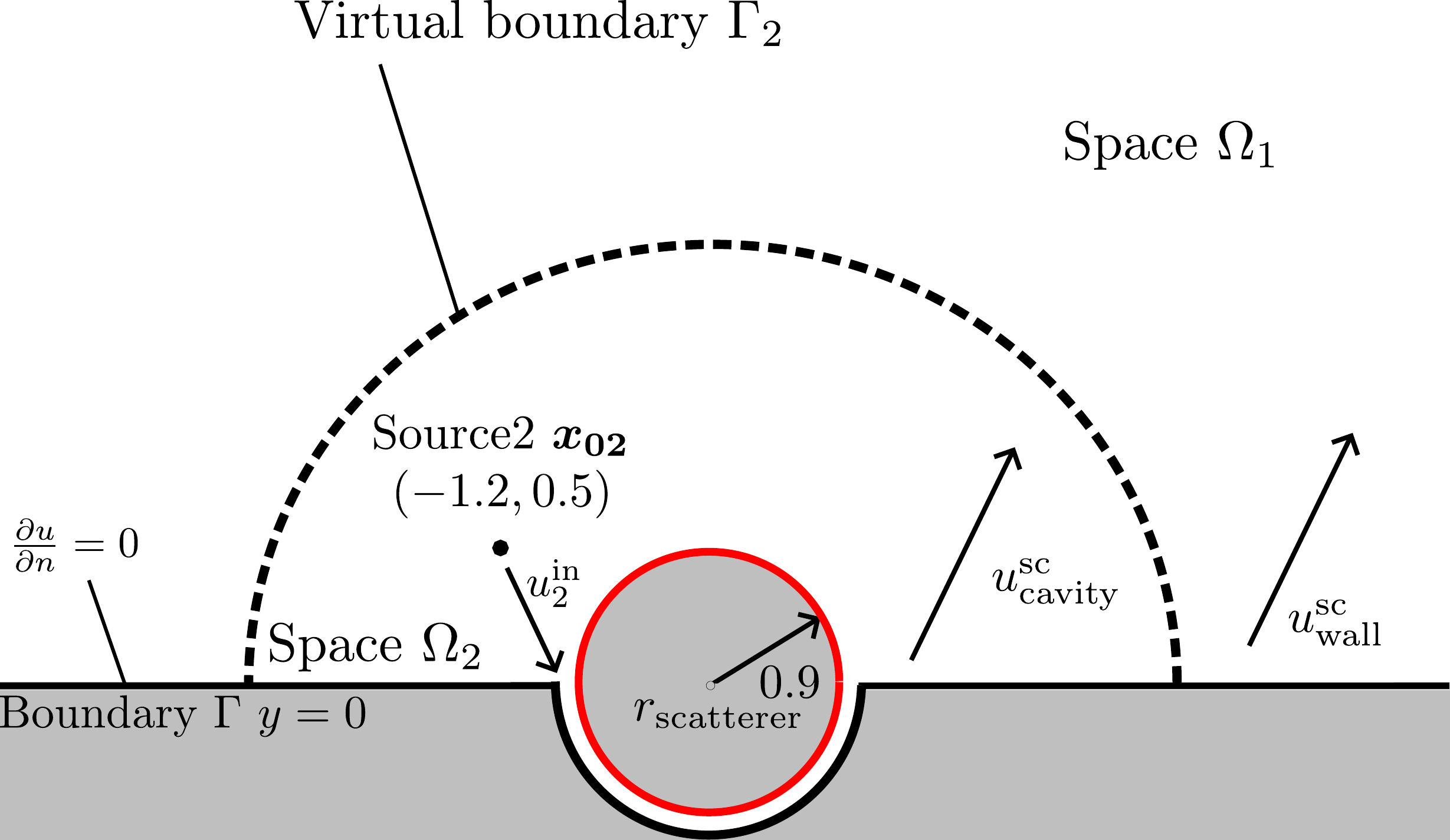}
 \caption{Problem settings for the case that a scatterer is arranged in a cavity.}\label{fig:5_setting4}
\end{figure}
\begin{figure}[h]
  \centering
  \includegraphics[scale=0.55]{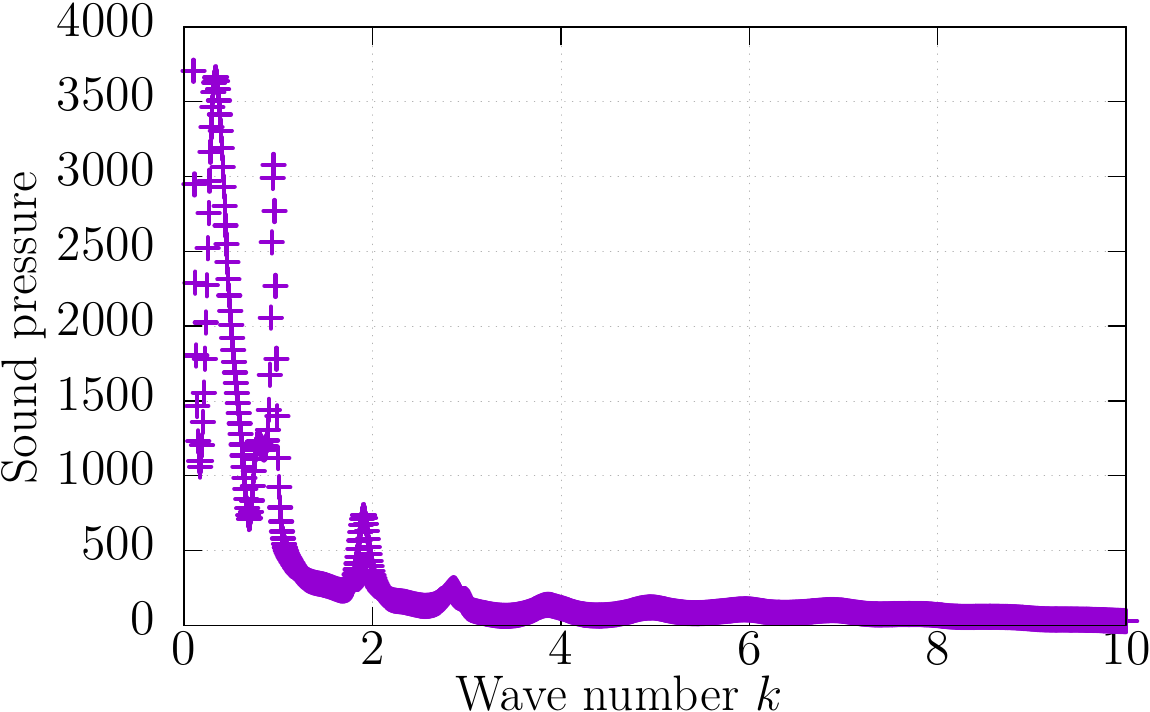}
  \caption{Sum of sound intensity at the right half of the distributed inner points.}\label{fig:5_suml2-2}
\end{figure}
\begin{figure}[h]
  \centering
  \includegraphics[scale=0.15]{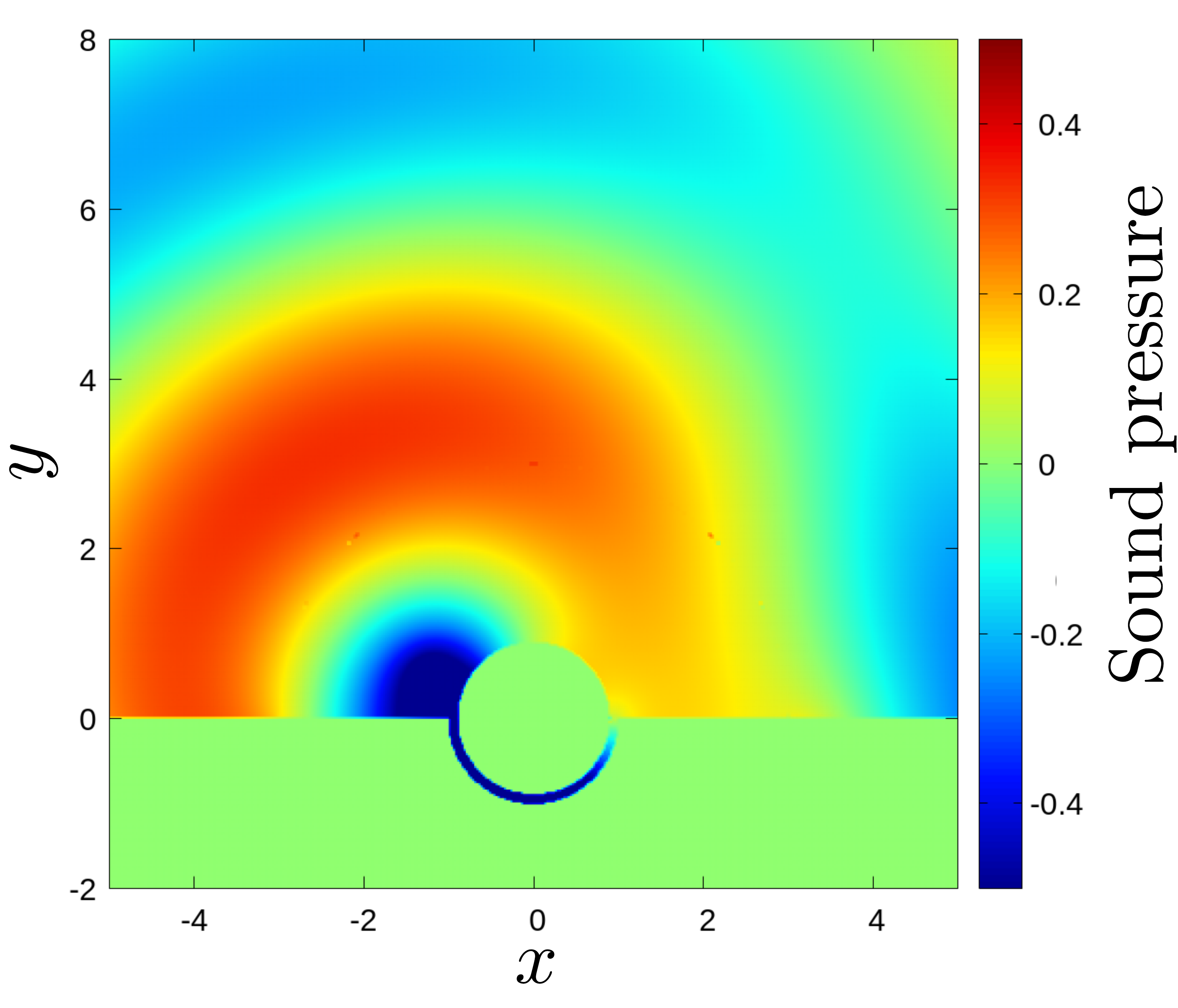}
  \caption{Sound pressure distribution with the incident wave number $k=0.78$}\label{fig:5_inhole}
\end{figure}

Figure \ref{fig:5_setting4} shows the setting for the second open resonator. We put a rigid circular scatterer centred at the origin and of the radius 0.9 in the domain used in the first validation. The domain is impinged by a point source located at $\bs{x_{02}}=(-1.2,0.5)$. The numerical settings are exactly the same as the ones in the previous example. By changing the wave number $k$ from $0.10$ to $10.00$ by $0.01$, we compute the total sound pressure at 160801 grid points in $-5.0 \le x \le 5.0$ and $-2.1 \le y \le 8.0$. Figure \ref{fig:5_suml2-2} shows the sum of the sound intensity at inner points which are located right half of the rectangular region.
In \figref{fig:5_suml2-2}, the graph again shows several sharp peaks. Figure \ref{fig:5_inhole} shows the sound pressure distribution $\Re[u(\bs{x})]$ for $\omega=0.78$. 
As indicated in \figref{fig:5_inhole}, the acoustic wave propagates through the region sandwiched by the cavity and the scatterer as if it works as a narrow waveguide at this excitation frequency. 

As the numerical demonstrations indicate, the proposed method allows us to solve acoustic scattering problems containing a cavity in a half-space that may work as an open resonator.

\section{Conclusion}
In this paper, we discussed a boundary element formulation for two-dimensional half-space and cavity scattering based on the hybrid integral representation that utilises both the layer potential and the Sommerfeld integral. We first modified the original formulation~\cite{lai2018new} to avoid the fictitious eigenvalue problem. The Burton-Miller formulation successfully ameliorates the accuracy of the proposed BEM. We also discussed appropriate settings for the parameter that splits the layer potential and the Sommerfeld integral. We then proposed a novel BEM for cavity scattering based on the hybrid BEM. The proposed formulation allows us to address complicated scattering that involves open resonances in half-space.

The straightforward future directions include the extension of the proposed method for three-dimensional half-space and cavities involving electromagnetic and elastic waves. In addition, the acceleration need also to be addressed, such as the fast Fourier transform (FFT), fast multipole method(FMM)~\cite{rokhlin1985rapid,greengard1987fast} and ${\cal H}$-matrix method~\cite{bebendorf2000approximation}, etc and/or exploit the block structure of the underlying matrix by using the Woodbury formula. As other possibilities, we can mention the topology and shape optimisations~\cite{bendsoe1988generating,isakari2017levelsetbased} in the half-space aiming to design a desired open resonator and related waveguide. The eigenvalue analysis~\cite{matsushima2022topology} for the integral equation is also of interest. They shall be discussed in the future elsewhere. 

\appendix 
\section{Integral representations for scattered fields in half-space}
 This appendix briefly reviews, for referential purposes, two integral representations for wave scattering in two-dimensional half-space. After stating the scattering problem of interest in Section \ref{sec:2.1}, we show the Sommerfeld and the hybrid integral representations~\cite{lai2018new} respectively in Sections \ref{sec:Sommer}, and \ref{sec:Hy_base}.
 \subsection{Settings}\label{sec:2.1}
\begin{figure}[h]
   \centering
   \includegraphics[scale=0.33]{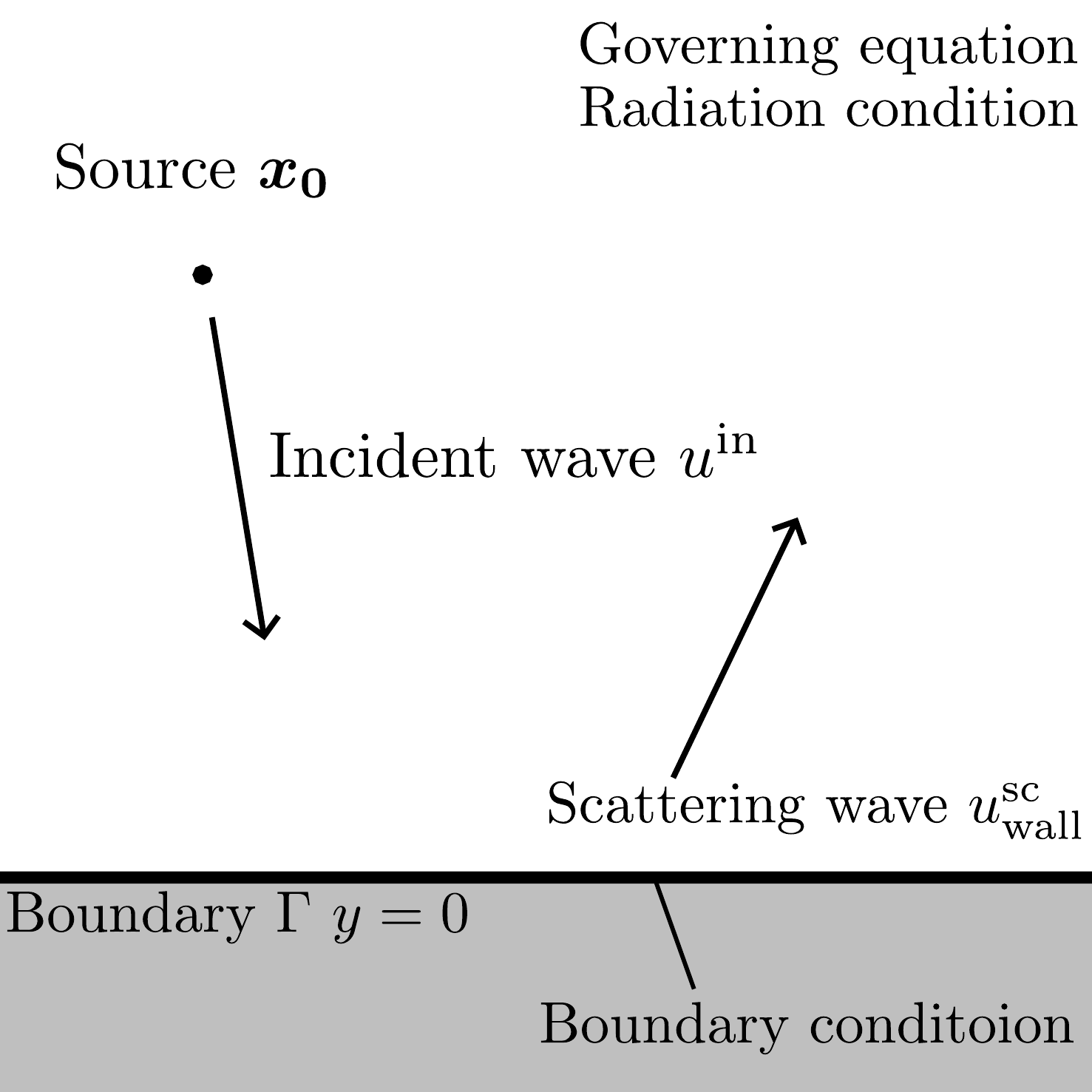}
   \caption{Acoustic scattering in a half-space.}\label{fig:setting1}
 \end{figure}

 In this article, we deal with harmonic acoustic scattering problem in the half-space $\mathbb{R}^{2}_{+}$~(\figref{fig:setting1}). The total field in this case solves the same BVP as \eqref{eq:Helmholtz}--\eqref{eq:radiation} except that the domain $\Omega$ is replaced by $\mathbb{R}^2_{+}$. 

\subsection{The Sommerfeld integral}\label{sec:Sommer}
It is well known that the scattered field for the problem in Section \ref{sec:2.1} is represented as
\begin{align}
  u^\mathrm{sc}(\boldsymbol{x})=\frac{1}{4 \pi} \int^{\infty}_{-\infty} \frac{e^{- \sqrt{\lambda^2 - k^2 }y}}{\sqrt{\lambda^2 -k^2}}e^{i \lambda x} \hat{\xi}(\lambda)  d\lambda, \label{eq:Sommer}
\end{align}
where $\hat{\xi}$ is a density function, and the square root is defined as
\begin{align}
 \sqrt{\lambda^2-k^2}=\begin{cases}
		       \sqrt{\lambda^2-k^2} & \text{if~}\lambda>k\\
		       -\mathrm{i}\sqrt{k^2-\lambda^2} & \text{otherwise}
		      \end{cases}.
\end{align}
In practical computation for \eqref{eq:Sommer}, the integral path is deformed as
\begin{align}
\lambda(t):=t-\mathrm{i}\frac{\tanh(t)}{a}, 
\label{eq:contour_deformation}
\end{align}
with a parameter $a>0$ to avoid the singularity at $\lambda=\pm k$. We used the symbol $\hat{\cdot}$ in \eqref{eq:Sommer} to emphasise that the density function lives in a Fourier space. $\lambda$ in \eqref{eq:Sommer} thus indicates the Fourier parameter. In this article, we define the Fourier transform of a univariate function $f(x)$ along the wall $y=0$ and its inverse respectively as follows: 
\begin{align}
  \hat{f}(\lambda):={\cal F}[f](\lambda) &:=\int^{\infty}_{-\infty}~f(x)~e^{-i \lambda x} dx, \label{eq:def_fourier} \\
  f(x):={\cal F}^{-1}[\hat{f}](x)&:=\frac{1}{2 \pi} \int^{\infty}_{-\infty}~\widehat{f}(\lambda)~e^{i \lambda x} d\lambda.
\end{align}
The unknown density function $\hat{\xi}$ can be determined by the boundary condition. If the integrand of \eqref{eq:Sommer} is rapidly decaying (that is the case e.g. when the source $\bs{x_0}$ is sufficiently away from the wall), we may replace the integral range $\mathbb{R}$ in \eqref{eq:Sommer} by a finite interval. In such a case, the density function can be obtained numerically, and accordingly, the scattered field can be computed.

\subsection{A solution based on hybrid integral representation.}\label{sec:Hy_base}
The Sommerfeld integral \eqref{eq:Sommer} is, unfortunately, not necessarily rapidly converging. In such a case, we may need a huge truncated integral interval, that prohibits us from utilising the Sommerfeld integral representation to solve the half-space problem. Lai et al~\cite{lai2018new} proposed a method to overcome this issue. In the proposed method, instead of the solely Sommerfeld integral representation, the sum of the Sommerfeld integral and the layer potential is adopted. This enables us to represent the scattered field only with the integrals over finite intervals. Roughly speaking, the layer potential part corresponds to the scattering from the wall close to the source, and the Sommerfeld integral to the rest.

We here review in some detail the scattering analysis for the problem in Section \ref{sec:2.1} with the hybrid representation. The scattered field $u^\mathrm{sc}$ in this case admits the following representation: 
\begin{align}
 u^{\mathrm{sc}}(\bs{x})=S_{\Gamma}[\sigma_W](\bs{x}) +F_{(-\infty, \infty)}[ \widehat{\xi_W}](\bs{x}), 
\label{eq:HBsol0}
\end{align}
where $S_{\Gamma}[\sigma_W]$ is the single layer potential on $\Gamma$ defined as
\begin{align}  
 S_{\Gamma}[\sigma_W](\boldsymbol{x}):=\int_{\Gamma}G_k(\boldsymbol{x},\boldsymbol{x'})\sigma_W(\boldsymbol{x'})ds(\boldsymbol{x'}), 
\label{eq:sec2_4_1}
\end{align}
with the fundamental solution 
\begin{align}
G_k(\bs{x},\bs{x'})=\frac{\mathrm{i}}{4}H_0^{(1)}(k|\bs{x}-\bs{x'}|)
\label{eq:fundamentalsol} 
\end{align}
of the 2D Helmholtz equation and the unknown density function $\sigma_W$. Also, $F_{(-\infty,\infty)}[\widehat{\xi_W}]$ in \eqref{eq:HBsol0} indicates the following Sommerfeld integral:
\begin{align}
 F_{(-\infty,\infty)}[\widehat{\xi_W}](\boldsymbol{x}):=\frac{1}{4\pi}\int^{\infty}_{-\infty}\frac{e^{- \sqrt{\lambda^2 - k^2 }y}}{\sqrt{\lambda^2 -k^2}}e^{i \lambda x} \widehat{\xi_W}(\lambda) d\lambda, 
\label{eq:sec2_4_2}
\end{align}
where $\widehat{\xi_W}$ is also the unknown density. If the density functions are appropriately chosen, the integral intervals $\Gamma$ and $\mathbb{R}$ in \eqref{eq:HBsol0} can respectively be replaced by $\Gamma_0:=\{\bs{x}\mid -M_0<x<M_0,~y=0\}$ and $I_0:=(-N_0, N_0)$, where $0<M_0$ and $N_0<\infty$ are positive and finite numbers, as shall be seen later. In other words, the integral representation \eqref{eq:HBsol0} can be reduced to 
\begin{align}
 u^{\mathrm{sc}}(\bs{x})=S_{\Gamma_0}[\sigma_W](\bs{x}) +F_{I_0}[ \widehat{\xi_W}](\bs{x}).
\label{eq:HBsol}
\end{align}
On the practical choice of the parameters $M_0$ and $N_0$, see Section \ref{sec:M0N0}. 

\begin{figure}[h]
 \centering
 \includegraphics[scale=0.5]{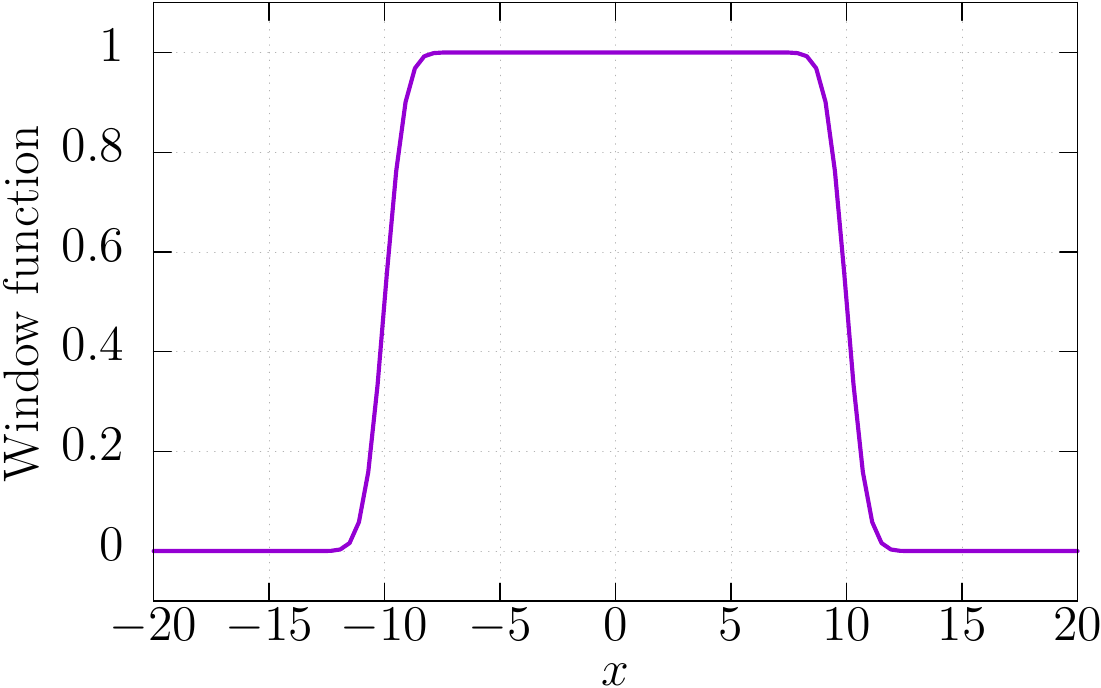}
 \caption{The window function in \eqref{eq:window} with $M_0=20$}\label{fig:window}
\end{figure}

The density functions $\sigma_W$ and $\widehat{\xi_W}$ are determined such that:
\begin{itemize}
\item[A.] the total field $u(\bs{x})=u^\mathrm{in}(\bs{x}) + u^\mathrm{sc}(\bs{x})$ satisfies the homogeneous Neumann boundary condition for $\bs{x}\in\Gamma$. 
\item[B.] both the integrands of the layer potential and the Sommerfeld integral \eqref{eq:HBsol} are locally supported. 
\end{itemize}
Here, following Lai at el \cite{lai2018new}, the density function $\sigma_W$ for the layer potential is chosen in light of the requirement B as
\begin{align}
\sigma_W(\bs{x})= W_{M_0}(x) \sigma(\bs{x}), \label{eq:sgmW}
\end{align}
with the window function $W_{M_0}$ defined as 
\begin{align}
 W_{M_0}(x):=\frac{1}{2} \left\{\mathrm{erf} \left( x+\frac{M_0}{2} \right)- \mathrm{erf} \left(x-\frac{M_0}{2} \right) \right\}, \label{eq:window}
\end{align}
with the error function $\mathrm{erf}(x)$. The window function in the case of $M_0=20$ is illustrated in \figref{fig:window}.
With the window function thus defined, $\mathrm{supp.}(\sigma_{W})\subset \Gamma_0$ is guaranteed. As $\sigma$ in \eqref{eq:sgmW}, we use a function solving the following integral equation:
\begin{align}
\frac{1}{2}\sigma(\bs{x})+D^*_{\Gamma_0}[\sigma](\bs{x}) + \frac{\partial u^\mathrm{in}}{\partial \bs{n}}(\bs{x}) = 0 
\label{eq:aug_int_eq}
\end{align}
for $\bs{x}\in\Gamma_0$, where $D^*_{\Gamma_0}[\sigma]$ indicates the normal derivative of the single layer potential defined as
\begin{align}
D^*_{\Gamma_0} [\sigma](\bs{x}) := \int_{\Gamma_0} \frac{\partial G(\bs{x},\bs{x'})}{\partial \bs{n}(\bs{x})} \sigma(\bs{x'}) ds(\bs{x'}).
\label{eq:dlayer}
\end{align}
The integral equation \eqref{eq:aug_int_eq} is inspired by the fact the single layer potential with density function $\sigma$ solving this equation satisfies the boundary condition on $\Gamma_0$. By using the obtained $\sigma_W$, the other density function $\widehat{\xi_W}$ is then determined by the requirement A as follows: 
\begin{align}
\frac{\partial u^\mathrm{in}}{\partial \bs{n}} + \frac{1}{2} \sigma_W + D^*_{\Gamma_0}[\sigma_W](\bs{x}) + \frac{1}{4\pi} \int_{-\infty}^{\infty} e^{\mathrm{i}\lambda x} \widehat{\xi_W} d\lambda = 0, \label{eq:global_equation}
\end{align}
which can easily be solved in the Fourier space as 
\begin{align}
 \widehat{\xi_W}=-2\widehat{\frac{\partial u^\mathrm{in}}{\partial \bs{n}}}-\widehat{\sigma_W}-2\widehat{D^{*}_{\Gamma_0}[\sigma_W]}.
 \label{eq:sol_global_equation}
\end{align}
Note that, owing to the vanishing layer potential $D^*_{\Gamma_0}[\sigma]$ in the current setting, both the density functions are obtained analytically as
\begin{align}
\sigma=-2\frac{\partial u^\mathrm{in}}{\partial \bs{n}}, \\
\widehat{\xi_W}={\cal F}[(1-W_{M_0})\sigma].
\end{align}
Since the window function is smooth, $\widehat{\xi_W}(\lambda)$ is rapidly decaying as $|\lambda|\rightarrow\infty$, from which all the requirements A and B are satisfied. It should be noted that the density functions can be chosen so that the integrals in \eqref{eq:HBsol} are defined over finite intervals even when some scatterers exist in the half-space and/or the boundary condition on $\Gamma$ is other than the homogeneous Neumann one. See the original paper~\cite{lai2018new} for further discussions. 

\paragraph{Acknowledgements}
At the early stage of this study, we communicated with Professor Jun Lai from Zhejiang University, who provided a careful guidance on the previous study. The authors also acknowledge that this work was supported by JSPS KAKENHI Grant Number 21K19764.


\begin{thebibliography}{10}
\expandafter\ifx\csname url\endcsname\relax
  \def\url#1{\texttt{#1}}\fi
\expandafter\ifx\csname urlprefix\endcsname\relax\def\urlprefix{URL }\fi
\expandafter\ifx\csname href\endcsname\relax
  \def\href#1#2{#2} \def\path#1{#1}\fi

\bibitem{sommerfeld1909uber}
A.~Sommerfeld, \href{https://doi.org/10.1002/andp.19093330402}{^^c3^^9cber die
  {Ausbreitung} der {Wellen} in der drahtlosen {Telegraphie}}, Annalen der
  Physik 333~(4) (1909) 665--736, publisher: John Wiley \& Sons, Ltd.
\newblock \href {https://doi.org/10.1002/andp.19093330402}
  {\path{doi:10.1002/andp.19093330402}}.

\bibitem{vanderpol1935theory}
B.~Van Der~Pol,
  \href{https://www.sciencedirect.com/science/article/pii/S0031891435901689}{Theory
  of the reflection of the light from a point source by a finitely conducting
  flat mirror, with an application to radiotelegraphy}, Physica 2~(1) (1935)
  843--853.
\newblock \href {https://doi.org/10.1016/S0031-8914(35)90168-9}
  {\path{doi:10.1016/S0031-8914(35)90168-9}}.

\bibitem{ochmann2004complex}
M.~Ochmann, \href{https://doi.org/10.1121/1.1819504}{The complex equivalent
  source method for sound propagation over an impedance plane}, The Journal of
  the Acoustical Society of America 116~(6) (2004) 3304--3311.
\newblock \href {https://doi.org/10.1121/1.1819504}
  {\path{doi:10.1121/1.1819504}}.

\bibitem{perez-arancibia2014highorder}
C.~P^^c3^^a9rez-Arancibia, O.~P. Bruno,
  \href{https://opg.optica.org/josaa/abstract.cfm?URI=josaa-31-8-1738}{High-order
  integral equation methods for problems of scattering by bumps and cavities on
  half-planes}, Journal of the Optical Society of America A 31~(8) (2014)
  1738--1746, publisher: Optica Publishing Group.
\newblock \href {https://doi.org/10.1364/JOSAA.31.001738}
  {\path{doi:10.1364/JOSAA.31.001738}}.

\bibitem{lai2018new}
J.~Lai, L.~Greengard, M.~O'Neil,
  \href{https://www.sciencedirect.com/science/article/pii/S1063520316300719}{A
  new hybrid integral representation for frequency domain scattering in layered
  media}, Applied and Computational Harmonic Analysis 45~(2) (2018) 359--378.
\newblock \href {https://doi.org/10.1016/j.acha.2016.10.005}
  {\path{doi:10.1016/j.acha.2016.10.005}}.

\bibitem{burton1971application}
A.~J. Burton, G.~F. Miller, \href{https://doi.org/10.1098/rspa.1971.0097}{The
  application of integral equation methods to the numerical solution of some
  exterior boundary-value problems}, Proceedings of the Royal Society of
  London. A. Mathematical and Physical Sciences 323~(1553) (1971) 201--210,
  publisher: Royal Society.
\newblock \href {https://doi.org/10.1098/rspa.1971.0097}
  {\path{doi:10.1098/rspa.1971.0097}}.

\bibitem{bao2014optimal}
G.~Bao, J.~Lai, \href{https://doi.org/10.1137/130905708}{Optimal {Shape}
  {Design} of a {Cavity} for {Radar} {Cross} {Section} {Reduction}}, SIAM
  Journal on Control and Optimization 52~(4) (2014) 2122--2140, publisher:
  Society for Industrial and Applied Mathematics.
\newblock \href {https://doi.org/10.1137/130905708}
  {\path{doi:10.1137/130905708}}.

\bibitem{bao2015radar}
G.~Bao, J.~Lai, {,Department of Mathematics, Zhejiang University, Hangzhou
  310027}, {,Courant Institute of Mathematical Sciences, New York University,
  New York, NY 10012},
  \href{http://aimsciences.org//article/doi/10.3934/dcdss.2015.8.419}{Radar
  cross section reduction of a cavity in the ground plane: {TE} polarization},
  Discrete \& Continuous Dynamical Systems - S 8~(3) (2015) 419--434.
\newblock \href {https://doi.org/10.3934/dcdss.2015.8.419}
  {\path{doi:10.3934/dcdss.2015.8.419}}.

\bibitem{ammari2000integral}
H.~Ammari, G.~Bao, A.~W. Wood,
  \href{https://doi.org/10.1002/1099-1476(200008)23:12<1057::AID-MMA151>3.0.CO;2-6}{An
  integral equation method for the electromagnetic scattering from cavities},
  Mathematical Methods in the Applied Sciences 23~(12) (2000) 1057--1072,
  publisher: John Wiley \& Sons, Ltd.
\newblock \href
  {https://doi.org/10.1002/1099-1476(200008)23:12<1057::AID-MMA151>3.0.CO;2-6}
  {\path{doi:10.1002/1099-1476(200008)23:12<1057::AID-MMA151>3.0.CO;2-6}}.

\bibitem{lai2014fast}
J.~Lai, S.~Ambikasaran, L.~F. Greengard,
  \href{https://doi.org/10.1137/140964904}{A {Fast} {Direct} {Solver} for
  {High} {Frequency} {Scattering} from a {Large} {Cavity} in {Two}
  {Dimensions}}, SIAM Journal on Scientific Computing 36~(6) (2014) B887--B903,
  publisher: Society for Industrial and Applied Mathematics.
\newblock \href {https://doi.org/10.1137/140964904}
  {\path{doi:10.1137/140964904}}.

\bibitem{kao2008maximization}
C.-Y. Kao, F.~Santosa,
  \href{https://www.sciencedirect.com/science/article/pii/S0165212507000923}{Maximization
  of the quality factor of an optical resonator}, Wave Motion 45~(4) (2008)
  412--427.
\newblock \href {https://doi.org/10.1016/j.wavemoti.2007.07.012}
  {\path{doi:10.1016/j.wavemoti.2007.07.012}}.

\bibitem{bao2016stability}
G.~Bao, K.~Yun, \href{https://doi.org/10.1007/s00205-015-0947-x}{Stability for
  the {Electromagnetic} {Scattering} from {Large} {Cavities}}, Archive for
  Rational Mechanics and Analysis 220~(3) (2016) 1003--1044.
\newblock \href {https://doi.org/10.1007/s00205-015-0947-x}
  {\path{doi:10.1007/s00205-015-0947-x}}.

\bibitem{rokhlin1985rapid}
V.~Rokhlin,
  \href{https://www.sciencedirect.com/science/article/pii/0021999185900026}{Rapid
  solution of integral equations of classical potential theory}, Journal of
  Computational Physics 60~(2) (1985) 187--207.
\newblock \href {https://doi.org/10.1016/0021-9991(85)90002-6}
  {\path{doi:10.1016/0021-9991(85)90002-6}}.

\bibitem{greengard1987fast}
L.~Greengard, V.~Rokhlin,
  \href{https://www.sciencedirect.com/science/article/pii/0021999187901409}{A
  fast algorithm for particle simulations}, Journal of Computational Physics
  73~(2) (1987) 325--348.
\newblock \href {https://doi.org/10.1016/0021-9991(87)90140-9}
  {\path{doi:10.1016/0021-9991(87)90140-9}}.

\bibitem{bebendorf2000approximation}
M.~Bebendorf, \href{https://doi.org/10.1007/PL00005410}{Approximation of
  boundary element matrices}, Numerische Mathematik 86~(4) (2000) 565--589.
\newblock \href {https://doi.org/10.1007/PL00005410}
  {\path{doi:10.1007/PL00005410}}.

\bibitem{bendsoe1988generating}
M.~P. Bends^^c3^^b8e, N.~Kikuchi,
  \href{https://www.sciencedirect.com/science/article/pii/0045782588900862}{Generating
  optimal topologies in structural design using a homogenization method},
  Computer Methods in Applied Mechanics and Engineering 71~(2) (1988) 197--224.
\newblock \href {https://doi.org/10.1016/0045-7825(88)90086-2}
  {\path{doi:10.1016/0045-7825(88)90086-2}}.

\bibitem{isakari2017levelsetbased}
H.~Isakari, T.~Kondo, T.~Takahashi, T.~Matsumoto,
  \href{https://www.sciencedirect.com/science/article/pii/S0045782516305187}{A
  level-set-based topology optimisation for acoustic^^e2^^80^^93elastic coupled
  problems with a fast {BEM}^^e2^^80^^93{FEM} solver}, Computer Methods in
  Applied Mechanics and Engineering 315 (2017) 501--521.
\newblock \href {https://doi.org/10.1016/j.cma.2016.11.006}
  {\path{doi:10.1016/j.cma.2016.11.006}}.

\bibitem{matsushima2022topology}
K.~Matsushima, H.~Isakari, T.~Takahashi, T.~Matsumoto,
  \href{https://www.sciencedirect.com/science/article/pii/S0165212522000634}{A
  topology optimization of open acoustic waveguides based on a scattering
  matrix method}, Wave Motion 113 (2022) 102987.
\newblock \href {https://doi.org/10.1016/j.wavemoti.2022.102987}
  {\path{doi:10.1016/j.wavemoti.2022.102987}}.

\end{thebibliography}

\end{document}